\documentclass[11pt]{amsart}
\usepackage{SSZ}
\usepackage[left=1.2in,top=1.5in,right=1.2in,bottom=1.5in]{geometry}

\renewcommand{\ss}{\mathord{\sslash}}

\begin{document}

\title{$K$-theory of finite monoid sets}
\author[Sarazola]{Maru Sarazola}
\address{School of Mathematics, University of Minnesota, Minneapolis MN, 55455}
\email{maru@umn.edu}
\author[Shapiro]{Brandon T. Shapiro}
\address{Department of Mathematics, Florida Gulf Coast University, Fort Meyers, FL, 33965}
\email{brandonshapiro@virginia.edu}
\author[Zakharevich]{Inna Zakharevich}
\address{Department of Mathematics, Cornell University, Ithaca, NY 14850}
\email{zakh@math.cornell.edu}

\begin{abstract}
In this paper we classify the algebraic $K$-theory of sets equipped with an action by a finite abelian monoid.  The main tool developed for this purpose is an Orbit-Stabilizer Theorem, which gives a decomposition of monoid sets analogous to the decomposition given by the Orbit-Stabilizer Theorem for $G$-sets.
\end{abstract}

\maketitle

\section*{Introduction}

There is an interesting observation, when working with abelian categories: that the category of finite sets bears a strange resemblance to an abelian category in far more situations than one might imagine.  It is not additive, but there is a notion of exact sequence (an inclusion and the inclusion of its complement) which gives a notion of kernel, cokernel, and image, and these notions exhibit similar dualities to those found in abelian categories.  Because of this, it is possible to define homology in finite sets \cite{ACGWhomology}, and this definition satisfies all of the same standard results such as the Snake Lemma as an abelian category.  It is tempting to conclude that there is a ``field with one element,'' generally called $\F_1$, such that finite sets are simply ``vector spaces over $\F_1$,'' and take this analogy seriously. This analogy can be taken surprisingly far.  There is a theory of ``schemes over $\F_1$'' which work analogously to ordinary schemes over a field, and the major theorems of algebraic geometry often hold in a similar manner.  For a detailed discussion of this, see for example \cite{deitmar05,F1schemes}.

The field $\F_1$ should be thought of as a field where ``you can't add, but you can multiply.''  So, for example, a finite set is a finite-dimensional $\F_1$-vector space, considered as the basis for the vector space.  If we remove all axioms of a vector space that refer to addition, and we can multiply only by the multiplicative identity, what remains are exactly sets.   On the other hand, a field extension over $\F_1$ is a set, containing $\F_1$, in which it is possible to multiply (but not to add), and in which every nonzero (i.e.\ ``not in $\F_1$''\footnote{The single element in $\F_1$ acts alternately like the multiplicative identity and like the zero element; it is a strange creature in this way.}) element has a multiplicative inverse.  In other words, it is an abelian group, with a disjoint basepoint added.  Understanding the vector spaces over this field extension is equivalent to understanding $G$-sets for a group $G$.  
The $K$-theory of a finite field extension of $\F_1$ is therefore simply the $K$-theory of $G$-sets for a finite abelian group $G$.  This has been known since \cite[Proposition 3.6]{segal_coh}:
\[K(G\mathbf{FinSets}) \simeq \bigvee_{H \leq G} \Sigma^\infty_+ B(G/H).\]
(This $K$-theory is also known without the assumption that $G$ is abelian, but in this paper we focus on the abelian case.)

For algebraic geometry, fields are not sufficient; it is important to understand algebras, as well.  What is an algebra over $\F_1$?  It is a set, with a basepoint, with an abelian multiplication and multiplicative unit (but no addition), with the basepoint acting as $0$.  In other words, it is a \emph{pointed abelian monoid}.  A module over such a monoid $A$ is a set $S$, equipped with a basepoint, with an action of $A$ for which the basepoint of $A$ acts as a zero; in other words, it is a (pointed) set with an action of $A$, analogous to a $G$-set.  The first main result of this paper is a formula for the $K$-theory of finite $A$-sets:
\begin{maintheorem} \label{main:mainthm}
    Let $A$ be a finite pointed abelian monoid.  The $K$-theory of $A$-sets splits as
    \[K(\ASets) \simeq \bigvee_{(R,F)\in OP_A} \Sigma^\infty_+ B(R\ss F).\]
\end{maintheorem}
Here $R$ and $F$ are unpointed submonoids of $A$ which make up an \emph{orbital pair}, and $R \ss F$ is their quotient. For the formal definition of a pointed abelian monoid, as well as associated terms, see  Subsection \ref{subsec.background}; for the definition of an orbital pair, see Definition~\ref{def:OP}; the stated result is Theorem \ref{mainthm}.  In particular, when $A$ is a group with a basepoint, this formula specializes to the previously-mentioned case of the $K$-theory of finite $G$-sets.

In order to properly understand the tools necessary for the proof of Theorem~\ref{main:mainthm}, let us consider the ingredients in the analysis of $K(G\FinSet)$.  The standard approach is the following.  First, we observe that any $G$-set splits into a disjoint union of $G$-orbits.  Secondly, we can classify $G$-orbits and show that any $G$-orbit is isomorphic to $G/H$ for some subgroup $H$.  Lastly, we show that there are no ``extensions'' of orbits, meaning that 
\[\iso G\FinSet \cong \prod_{H \leq G} \iso G\FinSet_{G/H},\]
where the terms of the product on the right are $G$-finite sets in which every orbit is $G$-equivariantly isomorphic to $G/H$.  Lastly, the Barratt--Priddy--Quillen--Segal Theorem is used to identify the terms on the right-hand side with $\Sigma^\infty_+ B(G/H)$.

When attempting to duplicate this analysis for monoid sets, every portion of it goes awry.   For a monoid set, the notion of an ``orbit'' is not cleanly defined, and in general it does not quite make sense.  Consider the following example, where the monoid $\N_+$ is acting on the set of nodes of the following graph by having $1$ take each node to the node it is pointing to, and having the basepoint $*$ take every point to $*$:
\begin{center}
    \begin{tikzpicture}
        \node (A) at (0,0) {$x$};
        \node (B) at (0,1) {$b$};
        \node (C) at (-1,0) {$c$};
        \draw[->] (B) -- (A);
        \draw[->] (C) -- (A);
        \draw (A) edge [loop right] (A);
        \node (pt) at (-1,1) {$*$};
        \draw (pt) edge [loop left] (pt);
    \end{tikzpicture}
\end{center}  If we na\"ively define the ``orbit'' of an element to be all of its images under $A=\N_+$, then the orbit of $b$ and the orbit of $c$ are not equal, but intersect at $x$---this is in direct contrast to the orbits of a group action, for which each pair of orbits is either disjoint or equal.  It may be tempting to simply declare an ``orbit'' to allow both ``forward'' and ``backward'' transversals of elements.  However, this definition breaks the \emph{next} portion of the analysis: that there are no ``extensions'' between orbits, and that therefore the category splits as a product of categories with particular orbit types.  The above picture gives a nontrivial extension of the $A$-subset $\{c,x\}$ to another $A$-set with a single ``orbit.'' 
\[
  \left(\begin{tikzpicture}[baseline={([yshift=-0.5ex]current bounding box.center)}]
        \node (A) at (0,0) {$x$};
        \node (C) at (-1,0) {$c$};
        \draw[->] (C) -- (A);
        \draw (A) edge [loop right] (A);
        \node (pt) at (-1,1) {$*$};
        \draw (pt) edge [loop left] (pt);
    \end{tikzpicture}\right)
    \rcofib
    \left(\begin{tikzpicture}[baseline={([yshift=-0.5ex]current bounding box.center)}]
                \node (A) at (0,0) {$x$};
        \node (B) at (0,1) {$b$};
        \node (C) at (-1,0) {$c$};
        \draw[->] (B) -- (A);
        \draw[->] (C) -- (A);
        \draw (A) edge [loop right] (A);
        \node (pt) at (-1,1) {$*$};
        \draw (pt) edge [loop left] (pt);
    \end{tikzpicture}\right)
    \rfib
    \left(\begin{tikzpicture}[baseline={([yshift=-0.5ex]current bounding box.center)}]
        \node (B) at (0,1) {$b$};
        \node (pt) at (-1,1) {$*$};
        \draw (pt) edge [loop left] (pt);
        \draw[->] (B) -- (pt);
    \end{tikzpicture}\right)
\]
This example also illustrates that even for very simple monoids $A$, classifying all $A$-``orbits'' would be extremely difficult.

The solution is to define a ``strong orbit'' to be a totally connected subset: a subset which can be transversed using only the action of $A$, but every step is ``reversible'' in the sense that if $a \cdot x = y$ for $a\in A$ then there exists $a'\in A$ with $a'\cdot y = x$. This definition was originally proposed by \cite{burnside}, and we explore it further in Section \ref{subsec:ambits}.  For strong orbits, there is also an Orbit-Stabilizer Theorem:
\begin{maintheorem}[Orbit-Stabilizer]
    If $A$ is an abelian monoid, then every $A$-set $S$ decomposes as a disjoint union of strong orbits.  Moreover, for every strong orbit $T \subseteq S$ there is an isomorphism of $A$-sets
    \[T_+ \cong (R_T \ss F_T)_+\]
    where $R_T \subseteq A$ is the submonoid of $T$-reversible elements, and $F_T\subseteq A$ is the stabilizer of any $x\in T$.
\end{maintheorem}
For a definition of $R_T$, $F_T$, and $\ss$, see Definitions~\ref{def:ss}, \ref{def:R_T} and \ref{def:F_T}; the stated result can be found in Proposition \ref{disjoint.union.of.ambits} and Theorem~\ref{ambit-stabilizer}. We also show in Theorem \ref{thm.uniqueness} that the submonoids $(R,F)$ that arise in this fashion are precisely the orbital pairs (Definition~\ref{def:OP}).

Finally, in Section \ref{functoriality} we study the functoriality of the computation of $K(\ASets)$ from Theorem \ref{main:mainthm} with respect to certain morphisms of monoids $A\to B$.   Section \ref{applications} contains some applications of our main results, several of which demonstrate that, in certain cases, the $K$-theory of monoid sets can be classified in terms of smaller monoids:

\begin{maintheorem} $\hbox{ }$
    \begin{enumerate}
        \item The $K$-theory of a monoid arising from a group is the same as the $K$-theory of the group.  More rigorously, if $G$ is an abelian group, then 
        \[K(G_+\text{-}\Sets) \simeq K(G\text{-}\Sets).\]
        \item Suppose that $A$ is a finite abelian monoid and $I \subseteq A$ is a nilpotent ideal.  Then 
        \[K(\ASets) \simeq K(A/I\text{-}\Sets).\]
        \item Suppose that $A$ is a finite local commutative algebra over a finite field $\mathbb{k}$. Then \[K(\ASets) \simeq K(A^\times\text{-}\Sets).\]
    \end{enumerate}
\end{maintheorem}
See Propositions~\ref{K(G)}, \ref{nilpotent.quotient}, and ~\ref{prop.comm.local.algebras} for more details.

\subsection*{Acknowledgments} The authors would like to thank Reid Barton, for help with Proposition~\ref{lem:injact}. The authors would also like to thank Christian Haesemeyer for many discussions, and in particular for asking the questions of Proposition \ref{transport} and Lemma \ref{nilpotent}. The authors are also grateful to Hiro Lee Tanaka for asking about the connection between our results and finite algebras, which led us to  Proposition \ref{prop.comm.local.algebras}.

The first author was partially supported by NSF grant DMS-2506116.  The last author was partially supported by NSF grant DMS-2405035.

\section{Strong orbits: definition and first results}

In this section we introduce the main characters: monoids, their actions, and the appropriate notion of orbit in this context. 

\subsection{Brief background on (pointed) monoids}\label{subsec.background}

We begin by setting some terminology and basic definitions, included here for the reader's convenience.

\begin{definition}\label{defn:ptdmonoid}
    A \emph{monoid} is a pointed set $A$ together with a binary operation $\cdot\colon A\times A\to A$ which is
    \begin{enumerate}
        \item associative; that is, $a\cdot (b\cdot c)=(a\cdot b)\cdot c$ for all $a,b,c\in A$;
        \item unital; that is, there exists an element $1\in A$ with $1\cdot a=a=a\cdot 1$ for all $a\in A$;
        \item such that $*\cdot a=*=a\cdot *$ for all $a\in A$, where $*$ denotes the basepoint of $A$.
    \end{enumerate}
    As usual, we say that a monoid $A$ is \emph{abelian} if $a\cdot b=b\cdot a$ for all $a,b\in A$.
\end{definition}

The monoids that concern this paper are always pointed, but their unpointed counterparts arise naturally in our results as well. For clarity, we include this definition below, and we also distinguish between the two possible types of substructures.

\begin{definition}
    An \emph{unpointed monoid} is a set $M$ together with a binary operation $\cdot\colon M\times M\to M$ which is
   associative and unital.
\end{definition}

\begin{definition}
    Let $A$ be a (pointed or unpointed) monoid. A subset $B\subseteq A$ is an \emph{unpointed submonoid} if  $1\in B$,  $B$ is closed under multiplication, and additionally $*\notin B$ in the case where $A$ is pointed. We say that $B$ is a \emph{pointed submonoid} of the pointed monoid $A$ if $*,1\in B$ and $B$ is closed under multiplication.
 \end{definition}

Unpointed submonoids give rise to a congruence relation on $A$, and hence to a notion of quotient; for the reader's convenience we recall this construction here.

\begin{definition} \label{def:ss}
    Let $F$ be an unpointed submonoid of a (pointed or unpointed) abelian monoid $A$.  Define the equivalence relation $\sim_F$ by $a\sim_F a'$ if there exists $x,x'\in F$ with $ax = a'x'$.  Write $A\ss F$ for the set of equivalence classes $A/\sim_F$.  When $F$ is clear from context we write $\sim$ instead of $\sim_F$.
\end{definition}

\begin{lemma}\label{lem.quotientA/F}
   If $A$ is abelian and $F\leq A$ is an unpointed submonoid, then $A\ss F$ is a well-defined abelian monoid, which is pointed or unpointed depending on whether $A$ is.
\end{lemma}

\begin{proof}
    We write $\sim$ for $\sim_F$ for brevity, and start by checking that $\sim$ is an equivalence relation. First, for any $a\in A$ we have that $a\sim a$ by setting $x = x' = 1\in F$.  Next, if $a \sim b$ then there exist $x,x'\in F$ with $ax = bx'$, and thus $bx' = ax$, so $b \sim a$. Lastly, suppose that $a \sim b$ and $b \sim c$.  Then there exist $x,x'\in F$ with $ax = bx'$ and $y,y'\in F$ with $by = cy'$.  But then \[axy = bx'y = x'by=x'cy'=cx'y';\] since $xy, x'y'\in F$, we have that $a \sim c$.  Thus $\sim_F$ is an equivalence relation, and $A\ss F$ is a well-defined set.

    To show that there is a well-defined multiplication, it suffices to show that if $a \sim b$ then $ac \sim bc$ for any $c\in A$.  But if $ax = bx'$ for some $x\in F$, then $acx = axc = bx'c=bcx'$, so $ac \sim bc$.  Thus $A\ss F$ inherits a well-defined multiplication from that of $A$, which is associative, unital, and abelian, by letting $[a]\cdot [b]\coloneqq [a\cdot b]$; we have just shown that this definition is independent of the choice of representatives in the equivalence classes.  %It also has a unit given by $[1]$, since if $1 \sim a$ then $r = ar'$ for some $r,r' \in R$.  Then for any $b\in A$, $[b][1] = [b]$, because $bar' = br \sim b$.
    In particular, if $A$ is pointed then so is $A\ss F$ with basepoint given by $[*]$, since by definition we have that $[*]\cdot [a]=[*\cdot a]=[*]$ for all $a\in A$.
\end{proof}

\begin{example}\label{groups.as.monoids}
   Any group $G$ gives rise to a (pointed) monoid by considering $G_+\coloneqq G\amalg \{*\}$, where the group multiplication is extended to $G_+$ by letting $*\cdot g=*=g\cdot *$ for all $g\in G_+$. If $F$ is a subgroup of an abelian group $G$, considered as an unpointed submonoid of $G_+$, then we have $g\sim_F g'$ if $xg = x'g'$, which is equivalent to saying $(x')^{-1}x = g'g^{-1}$---in other words, if $g$ and $g'$ are in the same coset of $F$. Then $A\ss F = (G/F)_+$; in particular, $F$ is not collapsed to the basepoint of $A$.  Thus the quotient is taken in an ``unbasepointed'' manner which may seem counterintuitive.
\end{example}

We now introduce our next main protagonists: sets equipped with a monoid action.

\begin{definition}
    Let $A$ be a monoid.  An \emph{$A$-set} is a pointed set $S$ together with a map $\cdot\colon  A \times S \rto S$ which is associative, unital, and pointed in the sense that
    \begin{enumerate}
    \item  $a \cdot (b\cdot x) = (ab) \cdot x$ for all $a,b\in A$ and $x\in S$;
        \item $1 \cdot x =x$ for all $x\in S$;
        \item $a \cdot * = *$ for all $a\in A$;
        \item $*\cdot x=*$ for all $x\in S$.
    \end{enumerate}
\end{definition}

\begin{definition}
    Let $A$ be a monoid. A \emph{morphism of $A$-sets} is a basepoint-preserving function $f\colon S\to S'$ between $A$-sets $S,S'$ which is $A$-equivariant, in the sense that $f(a\cdot x)=a\cdot f(x)$ for all $a\in A$ and all $x\in S$. 
\end{definition}

\begin{example}\label{ex:Gsets.as.Asets}
     If $G$ is a group, then any $G$-set $S$ gives rise to a $G_+$-set by considering $S_+\coloneqq S\amalg\{*\}$ and extending the action of $G$ on $S$ to an action of $G_+$ on $S_+$ by letting $g\cdot *=*$ for all $g\in G_+$ and $*\cdot x=*$ for all $x\in S_+$. 
\end{example}

\begin{example} \label{lem:multwelldef}
  The underlying pointed set $S=A$ of any monoid $A$ is always an $A$-set, with action given by (left) multiplication. Moreover, for any unpointed submonoid $F\leq A$ we have that $A \ss F$ is also an $A$-set with action defined as $a \cdot [b] = [a\cdot b]$, as we have already proved in Lemma \ref{lem.quotientA/F} that whenever $b \sim_F b'$ then for all $a\in A$, $ab \sim_F ab'$. 
\end{example}

\begin{example}\label{ex:A_d}
Consider the monoid \[A_d=\{*,1,t,t^2,\dots, t^{n-1} \mid t^n=t^d \text{ for some fixed } d\geq 1\}.\] This is a finite abelian monoid with the evident multiplication, which does not arise from a group. Moreover, this is not partially cancellative in the sense of \cite[Definition 1.2]{Hchuck}, as we have that $t^{n-d}t^d=1t^d\neq *$ but $t^{n-d}\neq 1$. %\maru{added this ex since our last meeting. I think this is a good example bc it's finite, abelian, but it has a weird cycle so none of the standard stuff applies}
\end{example}

%\begin{definition}
%        An $A$-set $S$ is \emph{finitely generated} if there exists a finite subset $S'\subseteq S$ such that for all $x\in S$ there exists $x'\in S'$ and $a\in A$ such that $x = a \cdot x'$.
%\end{definition}

%NOTE FOR US: since our main result is for finite monoids A, I think we might as well restrict to that since the start. Then, an A-set is finitely generated if and only if it's finite, so we can just deal with finite A-sets when needed and not introduce the finitely generated terminology. Similarly, I don't think it's worth introducing any Noetherianness conditions, and just deal with A finite. 

\subsection{Strong orbits}\label{subsec:ambits}

As explained in the introduction, the notion of \emph{orbit} is not well-behaved in the setting of $A$-sets, due to the possible lack of inverses in the monoid $A$. %Instead, we propose that the analogous notion is that of \emph{ambits}, whose definition requires the following. 
The meaningful notion to consider in this context is that of \emph{strong orbits}, introduced in \cite[Section 2.3.3]{burnside} and recalled below. Although \cite{burnside} works with unpointed monoids and partially-defined actions, his approach is equivalent to our pointed one and all results can be translated; see \cite[Section 1.7]{burnside} for a brief discussion.  The results in \cite{burnside} %introduced the notion of strong orbits, it focus on
deal with the very general case of monoids, instead of focusing on the special case of abelian monoids.  Here we exploit the strength of the abelianness assumption in order to produce a complete classification of strong orbit types. 
 
 \begin{definition}\label{defn.strongorbit}
 Let $A$ be a monoid, and $S$  be an $A$-set.   Define a relation $\sim$ on the elements of  $S$ by $x\sim y$ if and only if
there exist $a,b\in A$ such that $a\cdot x=y$ and $b\cdot y=x$. One can check that $\sim$ is an equivalence relation, and we call its
equivalence classes the \emph{strong orbits} of $S$.
\end{definition}

We find it especially helpful to provide an alternative characterization of strong orbits aside from the above one put forth in \cite{burnside}. Instead of describing strong orbits in terms of an equivalence relation, we characterize the strong orbit of an element $x\in S$ as the set obtained by acting on $x$ by a certain submonoid of $A$. We now present the required definitions.

\begin{definition}
    Let $A$ be an abelian monoid, $S$ an $A$-set, and let $x\in S$.  An element $a\in A$ is \emph{$x$-reversible} if there exists $b\in A$ such that $ba\cdot x = x$.    %An element $a\in A$ is \emph{$x$-eliminating} if $ax = *$. Let $E_x \subseteq A$ be the subset of $x$-eliminating elements of $A$.  
    An element $a\in A$ \emph{fixes $x$} if $a \cdot x = x$. 
    \end{definition}

\begin{notation}
    We denote by $R_x$  the subset of $A$ consisting of the $x$-reversible elements, and by $F_x$ the subset of $A$ consisting of the elements that fix $x$.  
\end{notation}

%The following is clear from the definition.

%\begin{lemma}
%    $E_x$ is an ideal, in the sense that if $a\in E_x$ then for all $b\in A$, $ba \in E_x$.   
%\end{lemma}

\begin{definition}\label{defn.ambits}
  Let $A$ be an abelian monoid, $S$ an $A$-set, and let $x\in S$.  The \emph{strong orbit of $x$} is the set $R_x \cdot x=\{a\cdot x \mid a\in R_x\}$.
\end{definition}

We now check that Definitions \ref{defn.strongorbit} and \ref{defn.ambits} are equivalent characterizations of the same notion.

\begin{lemma}
     Let $A$ be an abelian monoid, $S$ an $A$-set, and let $x\in S$. Then, the strong orbit containing $x$ as introduced in Definition \ref{defn.strongorbit} is precisely the set $R_x \cdot x$.
\end{lemma}
\begin{proof}
    Suppose that an element $y\in S$ is in the strong orbit of $S$ containing $x$ as introduced in Definition \ref{defn.strongorbit}; that is, there are elements $a,b\in A$ with $a\cdot x=y$ and $b\cdot y=x$. Then $ba\cdot x=b\cdot(a\cdot x)=b\cdot y=x$, which shows that $a\in R_x$ and thus $y\in R_x\cdot x$. Conversely, given $y=a\cdot x$ with $a\in R_x$, there exists some $b\in A$ such that $ba\cdot x=x$. Hence $b\cdot y=x$ and we conclude that $y$ is in the strong orbit containing $x$ in the sense of  Definition \ref{defn.strongorbit}.
\end{proof}

\begin{example}
    For any abelian monoid $A$ and $A$-set $S$, we have that $R_*=A$. Indeed, for any $a\in A$ we have $a\cdot *=*$, hence $1a\cdot *=*$ and $a$ is $*$-reversible. The strong orbit of $x=*$ is always the singleton set $\{*\}$. 
\end{example}

\begin{example}
    Let $G$ be an abelian group and $S$ be a $G$-set, and consider the $G_+$-set $S_+$ as in Example \ref{ex:Gsets.as.Asets}. Then for any $*\neq x\in S_+$ we have that $R_x=G$, as $g^{-1}g\cdot x=x$ for all $g\in G$. The strong orbit of a non-basepoint element $x\in S_+$ is thus $G\cdot x$; i.e.\  its orbit. 
\end{example}

\begin{example}\label{ex:A_d.ambits}
Let $A_d$ be the abelian monoid from Example \ref{ex:A_d}. Then for any $x\in\{1,t,\dots, t^{d-1}\}$ we have that $R_x=\{1\}$, but for any other  element $x\in\{t^d,\dots,t^{n-1}\}$ we see that $R_x=A_d\setminus \{*\}$, since this subset forms a cycle under the $A_d$-action. The strong orbit for an element $x$ is then either $\{x\}$ if $x\in\{1,t,\dots, t^{d-1}\}$, or $\{t^d,\dots,t^{n-1}\}$  otherwise. 
\end{example}

As illustrated by the examples above, strong orbits have some properties that mirror those of orbits in the classical setting; we establish some of these preliminary properties now.

\begin{proposition}\label{disjoint.union.of.ambits}
  Let $A$ be any abelian monoid, and $S$ be an $A$-set. Then $S$ decomposes as a disjoint union of strong orbits.  Moreover, if $x,y\in S$ are such that $R_x \cdot x = R_y \cdot y$ then $R_x = R_y$.
\end{proposition}
\begin{proof}
   % Write $S_x = R_x \cdot x$ for brevity. We know that $S = \bigcup_{x\in S} S_x$, since $1\in R_x$ for all $x$ and so $x\in S_x$.  Thus it suffices to check that for $x,y\in S$, either $S_x = S_y$ or $S_x \cap S_y = \emptyset $.

    %Suppose that $z\in S_x\cap S_y$.  Then $z = a\cdot x$, and since $a\in R_x$ that means that there exists $b\in A$ such that $b\cdot z = x$.  Similarly, there exist $a',b'\in R_y$ such that $a' \cdot y = z$ and $b' \cdot z = y$.  In particular, we see that $b'a\cdot x = y$, where $b'a\in R_x$ as we have that \[(ba')(b'a)\cdot x=ba'\cdot y=b\cdot z=ba\cdot x=x.\]  Thus $y\in S_x$, and hence $S_y \subseteq S_x$.  By symmetry, $S_x \subseteq S_y$, and we conclude that they are equal.  Thus two ambits are either equal or disjoint, as desired.

The fact that $S$ is a disjoint union of its strong orbits is immediate from Definition \ref{defn.strongorbit}, as strong orbits are equivalence classes under an equivalence relation. 

    Now suppose that $R_y \cdot y \subseteq R_x \cdot x$, and so there is some $a\in R_x$ with $a \cdot x = y$.  Suppose that $a'\in R_x$, so that there exists $b'$ with $b'a'\cdot x = x$; we claim that $a'\in R_y$.  Indeed, we have 
    \[b'a' \cdot y = b'a' \cdot (a\cdot x) = a \cdot (b'a'\cdot x) = a\cdot x = y.\]
    Thus $R_x \subseteq R_y$. By symmetry, we conclude that if $R_y \cdot y = R_x \cdot x$ then $R_x = R_y$ as desired.
\end{proof}

The result above makes it possible to set the following notation.

\begin{definition} \label{def:R_T}
    For a strong orbit $T$, write $R_T\defeq R_x$ for any $x\in T$.
\end{definition}

Note that $A$ acts transitively on any strong orbit, in the sense that for any elements $x,y$ in the same strong orbit $T$, we have that $T=R_x\cdot x$ and so there is some $a\in R_x\subseteq A$ such that $y=a\cdot x$; similarly, $T=R_y\cdot y$ and we get that $x=b\cdot y$ for some $b\in R_y\subseteq A$. Unlike the situation for $G$-sets, it is not necessarily the case that each strong orbit itself is an $A$-set, since strong orbit may not be closed under the action of $A$. For instance, this was already illustrated in Example \ref{ex:A_d.ambits}, where we have two distinct strong orbits $\{1\}$ and $\{t\}$ (if we assume $d>1$) and yet acting by $t$ maps the first strong orbit to the second. 

In spite of this, the action of $A$ still respects strong orbits in the following sense.

\begin{lemma}\label{a.T.contained.in.ambit}
    Let $T$ be a strong orbit of an $A$-set $S$, and fix any $a\in A$.  Then $a \cdot T$ is contained in a strong orbit of $S$.  In addition, for all $x\in S$ and all $a\in A$,  we have that $R_x \subseteq R_{a\cdot x}$.
\end{lemma}

\begin{proof}
    The first claim follows because if $b \cdot x = y$ then $b\cdot (a\cdot x) = a\cdot (b\cdot x) =a \cdot y$, so two elements $x,y$ in the same strong orbit are mapped into a single strong orbit.

    For the second statement, suppose that $b\in R_x$.  Then there exists $b'$ such that $b'b\cdot x = x$.  But then $b'b\cdot (a\cdot x) = a \cdot (b'b\cdot x) = a\cdot x$, so $b\in R_{a\cdot x}$; this shows that $R_x \subseteq R_{a\cdot x}$.
\end{proof}

\begin{example}
   The inclusion of strong orbits obtained in the above lemma need not be an equality. We can see this in Example \ref{ex:A_d.ambits} by letting $T=\{t^{d-1}\}$, $T'=\{t^d,\dots, t^{n-1}\}$, and $a=t$; then $a\cdot T\subsetneq T'$. 
\end{example}

The strong orbits which are closed under the $A$-action will play a central role in our study, and so we recall some terminology for them from \cite{burnside}, and present a  new technical result.

\begin{definition}
    Let $A$ be an abelian monoid and $S$ be an $A$-set. A strong orbit $T$ in $S$ is  \emph{invariant} if for all $x\in T$ and $a\in A$ we have that either $a\cdot x=*$ or $a\cdot x\in T$. 
\end{definition}

\begin{remark}
    As an immediate consequence of the definition we see that if $T$ is an invariant strong orbit, then $T_+\coloneqq T\cup \{*\}$ is an $A$-set, and so we have an $A$-equivariant inclusion $T_+\mrto S$. 
\end{remark}

\begin{lemma}\label{A/Rx.acts.trivially}
    Let $T$ be a non-trivial invariant strong orbit of an $A$-set $S$. Then, for any $x\in T$ we have \[A \backslash R_x = \{a\in A \mid a\cdot x = *\}.\]
\end{lemma}
\begin{proof}
The fact that the set on the right is included in $A\backslash R_x$ is immediate, as such elements may never be $x$-reversible. Conversely, if $a\in A$ is such that $a\cdot x\neq *$, since the strong orbit $T$ is invariant we must have $a\cdot x\in T=R_x\cdot x$, so there exists some $r\in R_x$ such that $a\cdot x=r\cdot x$. But then, if $r'\in R_x$ is an element such that $r'r\cdot x=x$, then this same element exhibits that $r'a\cdot x=x$, thus $a\in R_x$. 
\end{proof}

We conclude this subsection with a few useful technical results. 

\begin{definition}
    A submonoid $B \leq A$ is \emph{saturated} if whenever $ab\in B$ then $a\in B$ and $b\in B$. 
 %   An element $x$ is \emph{maximal} if for any $a\notin R_x$, $a\cdot x = *$.
\end{definition}

\begin{lemma} \label{lem:saturatedsub}
  For any abelian monoid $A$, any $A$-set  $S$, and element $*\neq x\in S$, the sets  $R_x$ and $F_x$ are unpointed submonoids of $A$, and $R_x$ is saturated.
\end{lemma}

\begin{proof}
   We have that $1\in R_x$, because $1 \cdot x = x$. To see that $R_x$ is closed under multiplication, let $a,a'\in R_x$; then there exist $b,b'\in A$ witnessing that $a,a'$ are $x$-reversible (respectively).  Then
    \[(bb')(aa')\cdot x = (ba)(b'a')\cdot x  = x,\]
    so $aa'\in R_x$. Moreover, we see that $*\notin R_x$, as otherwise we would have $x=b*\cdot x=*\cdot x=*$, which contradicts the assumption that $x\neq *$. This shows that $R_x$ is an unbased submonoid of $A$. It is straightforward to check that $F_x$ is an unpointed submonoid as well.

    To see that $R_x$ is saturated, suppose that $ab\in R_x$. Thus there exists $c\in A$ such that $cab\cdot x = x$.  But then both $a$ and $b$ are $x$-reversible, so $a,b\in R_x$, as desired.
\end{proof}

\begin{lemma} \label{lem:fixinv}
  Let $S$ be an $A$-set, and fix elements $x\in S$ and $a\in A$. Then $F_x \subseteq F_{a \cdot x}$, and moreover, these sets are equal whenever $a\in R_x$.
\end{lemma}

\begin{proof}
    Suppose that $b\in F_x$.  Then $b\cdot (a\cdot x) = a\cdot (b\cdot x) = a \cdot x$, so $F_x \subseteq F_{a \cdot x}$.  If in addition we have that  $a\in R_x$, then there exists $a'\in R_x = R_{r\cdot x}$ such that $a'a \cdot x = x$.  But then any element $b\in F_{a\cdot x}$ satisfies \[b\cdot x=b\cdot (a'a\cdot x)=a'\cdot(b\cdot (a\cdot x))=a'\cdot (a\cdot x)=x\] which shows that $F_{a\cdot x}\subseteq F_x$. 
\end{proof}

As a consequence, we may set the following notation.

\begin{definition} \label{def:F_T}
    For an $A$-strong orbit $T$, write $F_T\defeq F_x$ for any $x\in T$.
\end{definition}

\section{Classification of invariant strong orbits}\label{section:classifying.ambits}

In the classification of $G$-sets for a group $G$, and subsequently in the analysis of the $K$-theory of $G$-sets, it is important to classify $G$-orbits: these are the transitive $G$-sets $X$, which one can prove are (up to isomorphism) of the form $G/\Stab(x)$ where $\Stab(x)\leq G$ denotes the \emph{stabilizer} of any element $x\in X$.   The situation for $A$-sets is more subtle, partly due to the fact that strong orbits need not be closed under the $A$-action, which makes their characterization inherently more complicated. In this section we give a classification of invariant strong orbits, reminiscent of the one for $G$-orbits.

\subsection{The Orbit-Stabilizer Theorem}

Given an $A$-set $S$ and a non-trivial invariant strong orbit $T\subseteq S$, the role of the stabilizer of an element $x\in T$ is played by the unpointed submonoid $F_x=F_T\subseteq A$, which is independent of the choice of element $x$ by Lemma \ref{lem:fixinv}. While it is possible to consider the quotient $A\ss F_T$ as an $A$-set in a natural way (see Example \ref{lem:multwelldef}), this will generally not be isomorphic to $T_+$. 

However, this is not the analogy we would expect. In the case for $G$-sets, the isomorphism is between a $G$-orbit $X=G\cdot x$ and $G/\Stab(x)$; hence, in the case of $A$-sets, we would hope for an isomorphism between an invariant strong orbit $T=R_x\cdot x$ and $R_x\ss F_x$. This is the main result of this section.

We start by detailing the $A$-set structure on the pointed set $(R_x\ss F_x)_+$.

\begin{proposition}\label{R/TasAset}
     Let $A$ be an abelian monoid, $S$ be an $A$-set, and $T\subseteq S$ be a non-trivial invariant strong orbit. Then the pointed set $(R_T\ss F_T)_+$ admits an $A$-action given by \[a \cdot [r] = \begin{cases}
            [ar] \caseif a \in R_T, \\
            [*] \caseif a \notin R_T.
    \end{cases}\] Moreover, this $A$-set has a single (and thus necessarily invariant) non-trivial strong orbit.
\end{proposition}
\begin{proof}
  By definition, we immediately have that $a\cdot [*]=[*]$ for all $a\in A$; that $1\cdot [r]=[1\cdot r]=[r]$ for all $[r]$ since $1\in R_T$; and that $*\cdot [r]=[*]$ for all $[r]$ since $*\notin R_T$. Hence, to check that this is a well-defined action it suffices to show that $ab \cdot [r] = a\cdot (b\cdot [r])$ for all $[r]$ and all $a,b\in A$. This clearly holds if both $a$ and $b$ are in $R_T$. If at least one of $a$ and $b$ is not in $R_T$, then $ab$ is not in $R_T$ (since $R_T$ is saturated by Lemma \ref{lem:saturatedsub}), so associativity holds as well. We conclude that we have a well-defined $A$-set.  
  
  To see that this $A$-set has a single non-trivial strong orbit, consider a non-basepoint element $[r]\in R_T\ss F_T$; we want to show that $[r]$ is in the same strong orbit as $[1]$. We know that $r\in R_T$, and so if $x\in T$ is any element we have that $r'r\cdot x=x$ for some $r'\in R_T$. In particular, this means that $r'r\in F_T$ and thus $r'r\cdot [1]=[r'r]=[1]$. Hence $r$ is $[1]$-reversible, and since $[r]=r\cdot [1]$, this shows that $[r]$ is in the same strong orbit as $[1]$, as desired.
\end{proof}

Our main result is then the following.

\begin{theorem}[Orbit-Stabilizer Theorem]\label{ambit-stabilizer}
    Let $A$ be an abelian monoid and $S$ be an $A$-set. Then, for every  non-trivial invariant strong orbit $T\subseteq S$ there is an isomorphism of $A$-sets $$T_+\cong (R_T\ss F_T)_+.$$
\end{theorem}
\begin{proof}
    We endow the quotient $(R_T\ss F_T)_+$ with an $A$-action as in Proposition \ref{R/TasAset}, and define maps in both directions as follows. Fix $x\in T=R_T\cdot x$, and let $f\colon T_+\to (R_T\ss F_T)_+$ be the pointed function given by $f(r\cdot x)=[r]$, and $g\colon (R_T\ss F_T)_+\to T_+$ be the pointed function given by $g([r])=r\cdot x$. We start by showing that both $f$ and $g$ are well-defined.

    To prove that $f$ is well-defined, suppose that $r\cdot x = s\cdot x$ for some $r,s\in R_T$.  We claim that this implies that $r\sim_{F_T} s$. Since $r\in R_T$, there is some $r'\in R_T$ such that \[x=r'r\cdot x = r'\cdot (r\cdot x)=r'\cdot (s\cdot x)=r's\cdot x;\] hence $rr',sr'\in F_T$.  In particular, $rr',sr' \sim_{F_T} 1$ and thus $s \sim_{F_T} rr's \sim_{F_T} r$, as desired.  
    
    Similarly, to show that $g$ is well-defined, suppose that $[r]=[s]$; we must check that $r\cdot x=s\cdot x$. Since $r\sim_{F_T} s$, there exist $r',s'\in F_T$ with $rr'=ss'$. Then \[r\cdot x=r\cdot (r'\cdot x)=rr'\cdot x=ss'\cdot x=s\cdot(s'\cdot x)=s\cdot x,\] as desired.
    
    Since the functions $f$ and $g$ are clearly inverses of each other, it remains to show that one of them is $A$-equivariant, as this will imply that the other one is as well.  We show that $f$ is $A$-equivariant. Given $a\in A$ and $r\cdot x\in T$, we have that
 \[a\cdot f(r\cdot x)=a \cdot [r] = \begin{cases}
            [ar] \caseif a \in R_T, \\
            [*] \caseif a \notin R_T
    \end{cases}= \begin{cases}
            f(ar\cdot x) \caseif a \in R_T, \\
            f(*)=f(r\cdot *)=f(r\cdot(a\cdot x)) \caseif a \notin R_T
    \end{cases}\]
    where in the very last equality we have used the fact that any element $a\in A\backslash R_x$ must be such that $a\cdot x=*$, by Lemma \ref{A/Rx.acts.trivially}. We see that $a\cdot f(r\cdot x)=f(a\cdot(r\cdot x))$ for all $a\in A$, which concludes the proof.
\end{proof}

\begin{remark}
    The above result is of course also true for the trivial invariant strong orbit $T=\{*\}$, except in this case, we have that $R_*=A=F_*$, and so the role of the quotient is played simply by $\{*\}$. 
\end{remark}

\begin{remark}%\brandon{added this}\maru{looks good to me!}
    For a strong orbit $T$ which is not invariant, the proof of Theorem~\ref{ambit-stabilizer} shows that the isomorphism $T_+\cong (R_T\ss F_T)_+$ still holds, except that the $A$-set structure on $T_+$ is no longer induced by the inclusion of sets $T_+ \subseteq S$. Instead, for $x \in T$ and $a \notin R_T$ we set $a \cdot x = *$. Equivalently, the isomorphism $T_+\cong (R_T\ss F_T)_+$ can be regarded as one of $R_T$-sets rather than $A$-sets without modifying the $R_T$-action on $T_+$.
\end{remark}

\begin{remark}
    In \cite[Proposition 3.1]{burnside}, the author gives a different characterization of strong orbits. They show that any strong orbit $T$ is isomorphic to a quotient of the form $R_e/\equiv$, where $e\in A$ is a minimal idempotent that fixes some element $x$ in the strong orbit, $R_e$ is the strong orbit of $e$ when we consider $A$ as an $A$-set with multiplication, and $\equiv$ is the congruence relation on $R_e$ defined by $a\equiv b$ if and only if $a\cdot x=b\cdot x$. Although their quotient $R_e/\equiv$ is, of course, isomorphic to our quotient $R_T\ss F_T$ of Theorem \ref{ambit-stabilizer} (as they are both isomorphic to the given strong orbit $T$), in general we have that $R_e\subsetneq R_T$, so these results provide two alternative characterizations of strong orbits. Their characterization has the benefit that the set $R_e$ is smaller than our $R_T$; on the other hand, our characterization does not require one to find an appropriate idempotent $e$ nor its orbit, and $R_T$ can be directly extracted for instance from a graphical representation of the $A$-set in question. A further advantage of our characterization will appear in Theorem \ref{thm.uniqueness}, where we describe all the possible pairs $(R_T, F_T)$ that give rise to invariant strong orbits in terms of submonoids of $A$, without the need to consider any $A$-sets.
\end{remark}

In the case of $G$-sets, the fact that every orbit is isomorphic to a $G$-set of the form $G/H$ for some subgroup $H\leq G$ has an additional immediate consequence: when $G$ is abelian every orbit has a group structure in a natural way. We recover an analogous result for $A$-sets.

\begin{proposition} \label{lem:quotisgp}
    Let $A$ be an abelian monoid and $S$ be an $A$-set. Then for every invariant strong orbit $T$, we have that  $R_T\ss F_T$ is an abelian group. 
    %Suppose that $A\ss F$ is a pointed ambit.  Then $(A \ss F) \backslash \{*\}$ is an abelian group. 
\end{proposition}

\begin{proof}
%It suffices to show that $R_T\ss F_T$ is an abelian group, with unit $[1]$ and with multiplication given by $[r] [s]=[rs]$.
%The fact that this multiplication is well-defined is a consequence of Lemma \ref{lem.quotientA/F}, together with the fact that $R_T$ is closed under multiplication by Lemma \ref{lem:saturatedsub}. Hence, this multiplication is associative, unital, and abelian (as it is inherited from $A$), and it only remains to show that every element $[r]\in R_T\ss F_T$ has an inverse.
It suffices to show that every element $[r]\in R_T\ss F_T$ has an inverse.
   Fix some element $x\in T$. Since $r\in R_T=R_x$, there is some $r'\in R_T$ such that $r'r\cdot x=x$; thus $r'r\in F_T=F_x$ and we have that $[1]=[r'r]=[r'][r]$. Since $A$ is abelian, this shows that $[r']=[r]^{-1}$, as desired. 
\end{proof}

%\begin{definition}
 %   An $A$-set $S$ \emph{has a single $A$-ambit} if for any two nonbasepoint elements $x,y\in S$ there exists an element $a\in A$ such that $a \cdot x = y$.  When $A$ is clear from context we simply say that \emph{$S$ has a single ambit}.
%\end{definition}

% \begin{definition}
%     An $A$-set $S$ is a \emph{pointed ambit} (or \emph{pointed $A$-ambit}, if $A$ is not clear from context) if $S = T_+$ for $T$ an $A$-ambit, and if for all nonbasepoint $a\in A$ and $x\in T$, $a\cdot t \neq *$. 
% \end{definition}

\subsection{A computation of all invariant strong orbit types} 
%We now investigate the extent to which the characterization of closed ambits in terms of quotients of the form $R_T\ss F_T$ given in Theorem \ref{ambit-stabilizer} is unique. Evidently, g
 %The goal of this subsection is to show that it is not just the isomorphism type of the quotient that is determined, but in fact, there is a unique pair of submonoids $(R,F)$ of $A$ which realizes this quotient. 
In this subsection, we isolate the precise conditions that a pair of submonoids $(R,F)$ of $A$ must satisfy in order for their quotient $R\ss F$ to model an invariant strong orbit type. This allows us to compute all possible invariant strong orbit types by analyzing the (finite!) monoid $A$, without the need to study all the possible $A$-sets. Moreover, given isomorphic invariant strong orbits $T\cong T'$, the Orbit-Stabilizer theorem gives an isomorphism $R_T\ss F_T\cong R_{T'}\ss F_{T'}$. We show that it is not just the isomorphism type of the quotient that is determined, but in fact, there is a unique pair of submonoids $(R,F)$ with the relevant conditions which realizes this quotient.

We start with the following useful lemma.

\begin{lemma} \label{lem:ambittopair}
  Let $A$ be an abelian monoid and $f:S \rto S'$ be an injection of $A$-sets.  Then $R_x = R_{f(x)}$ and $F_x = F_{f(x)}$ for all $x\in S$.  %In particular, the isomorphism type of an $A$-set with a single $A$-ambit uniquely determines a pair $(R,F)$. 
\end{lemma}

\begin{proof}
    For any $r\in R_x$, there exists some $r'\in R_x$ such that $r'r\cdot x = x$.  Then $r'r \cdot f(x) = f(r'r\cdot x) = f(x)$, so $r\in R_{f(x)}$; thus $R_x \subseteq R_{f(x)}$.  Similarly, if $r \cdot x = x$ then $r \cdot f(x) = f(x)$, so $F_x \subseteq F_{f(x)}$.

    On the other hand, suppose that $r\in R_{f(x)}$.  Then there exists $r'\in R_{f(x)}$ such that $r'r \cdot f(x) = x$.  Then $f(r'r\cdot x) = f(x)$, and since $f$ is injective it follows that $r'r \cdot x = x$, and thus $r\in R_{x}$; thus $R_{f(x)} \subseteq R_x$.  Similarly, $F_{f(x)} \subseteq F_x$.  
\end{proof}

Combining our results thus far, we can formulate the following result, which serves both to describe the uniqueness of the expression obtained in Theorem \ref{ambit-stabilizer}, and to characterize the pairs of submonoids $F\leq R\leq A$ whose quotient $R\ss F$ may appear as an invariant strong orbit.  Before stating the theorem we introduce some terminology.

\begin{definition} \label{def:OP}
Let $A$ be an abelian monoid.
    A \emph{orbital pair in $A$} is a pair $(R,F)$ of unpointed submonoids of $A$ satisfying the following properties:
    \begin{enumerate}
        \item $F \leq R$.
        \item $R$ is saturated.
        \item  $R \ss F$ is a group.
        \item  If $a\in A$ and $f\in F$ such that $af\in F$ then $a\in F$.
    \end{enumerate}
    When $A$ is clear from context we simply say that $(R,F)$ is a orbital pair.  Denote the set of orbital pairs in $A$ by $OP_A$.
\end{definition}

\begin{theorem}\label{thm.uniqueness}
   Let $A$ be an abelian monoid. There is a bijection between
    \begin{enumerate}
        \item The set of isomorphism classes of $A$-sets with a single non-trivial strong orbit
        \item The set $OP_A$ of orbital pairs in $A$. %pairs $(R,F)$ of unpointed submonoids $F\leq R\leq A$ such that $R$ is saturated, $R \ss F$ is a group, and such that for any $a\in A$ and $f\in F$, if $af \in F$  then $a \in F$.
    \end{enumerate}
\end{theorem}

\begin{proof}
    Suppose that $S$ is an $A$-set with a single non-trivial (thus necessarily invariant) strong orbit $T$.  For any element $x\in T$, we know by Lemma \ref{lem:saturatedsub} that $F_x\subseteq R_x \subseteq A$ are unpointed submonoids, and that $R_x$ is saturated. Moreover,  $R_x\ss F_x$ is a group by Lemma~\ref{lem:quotisgp}.  The submonoid $F_x$ has the further property that if $af\in F_x$ for some $a\in A$ and $f\in F_x$, then $a\in F_x$, since we have that $a\cdot x=a\cdot (f\cdot x)=af\cdot x=x$.    The assignment $S=T_+\mapsto (R_T, F_T)$ is well-defined within the isomorphism class of $S$, because if $T_+\cong T'_+$ then $(R_T,F_T) = (R_{T'},F_{T'})$ by Lemma~\ref{lem:ambittopair}.  This gives a function from (1) to (2).

Starting from a pair $(R,F)$ as in (2), we can define an $A$-action on the set $(R\ss F)_+$ as in Proposition \ref{R/TasAset}, for which $R\ss F$ is a non-trivial invariant strong orbit. Indeed, a detailed study of the proof of Proposition \ref{R/TasAset} reveals that what we need for it to hold is for $R$ to be saturated, and that for any $[r]\in R\ss F$ there is some $[r']\in R\ss F$ with $[r'r]=[1]$ (guaranteed here since $R\ss F$ is a group). We thus get a function from (2) to (1) given by $(R,T)\mapsto (R\ss T)_+$. 

The composite $T_+\mapsto (R_T, F_T)\mapsto (R_T\ss F_T)_+$ is the identity on isomorphism classes, by Theorem \ref{ambit-stabilizer}. To prove that the other composite $(R,F)\mapsto (R\ss F)_+\mapsto (R_{R\ss F}, F_{R\ss F})$ is the identity, 
 we must check that for any $[r]\in R\ss F$ we have that $R_{[r]} = R$ and $F_{[r]} = F$.  Suppose that $a\in A$ is $[r]$-reversible.  In particular, we must have $a \cdot [r] \neq [*]$, but by definition of the action of $A$ on $(R\ss F)_+$, this implies that $a\in R$, so $R_{[r]} \subseteq R$.  On the other hand, suppose that $a\in R$.  Since $R\ss F$ is a group, there exists $b\in R$ such that $[ba]=[1]$.  But then $ba\cdot [r] = [bar] = [ba][r]=[r]$, so $a\in R_{[r]}$.  Thus $R \subseteq R_{[r]}$, and they are equal, as desired.  Now suppose that $a\in F$.  Then $a \cdot [r] = [ar] = [a][r]=[1][r]=[r]$, so $F \subseteq F_{[r]}$.  On the other hand, suppose that $a\in F_{[r]}$, and so $a\cdot [r] = [r]$.  Since $R \ss F$ is a group, there exists $q\in R$ such that $[q][r]=[1]$.  Then 
 $[a]=a\cdot [1]=a\cdot [qr]=q\cdot (a\cdot [r])=q\cdot [r]=[1]$, or in other words $a \sim_F 1$.
 This is exactly the statement that there exists $f\in F$ such that $af\in F$, so by assumption this implies that $a\in F$.  Thus $F_{[r]} \subseteq F$, and they are equal, as desired.  
\end{proof}

\subsection{Automorphisms of invariant strong orbits}

When $G$ is an abelian group, one can show that the set of $G$-automorphisms of any transitive $G$-set $G/H$ is in bijection with the set of cosets $G/H$ itself. The goal of this subsection is to prove the corresponding result for invariant strong orbits. We note that \cite[Theorem 3.2]{burnside} also contains a computation of the group $\Aut(R_e/\equiv)$, which in the case when the monoid is abelian, can be translated in terms of a quotient of $R_e$. Instead of leveraging this result and connecting it to our description of strong orbits from Theorems \ref{ambit-stabilizer} and \ref{thm.uniqueness}, we choose to include a computation here using different methods, both to illustrate how the study can be simplified when the monoid is abelian, and to give the reader a self-contained account that does not require us to introduce the notation of \cite{burnside}.

We start with the following result, which gives us a way to construct bijections from a strong orbit to itself. 

\begin{proposition} \label{lem:injact}
    Let $A$ be an abelian monoid, $S$ be an $A$-set, and $T\subseteq S$ be any strong orbit.  Then for any element $r\in R_T$, the function $r\cdot: T \rto T$ is a bijection.  
\end{proposition}

\begin{proof}
First note that this function has the indicated codomain, as any $r\in R_T$ must be such that $r\cdot x\neq *$ for any $x\in T$ by Lemma \ref{A/Rx.acts.trivially}; thus $r\cdot x\in T$. 

    To check injectivity, suppose that $x,y\in T$ are such that $r \cdot x = r \cdot y$.  Let $r'\in R_T$ be such that $r'\cdot (r \cdot x) = x$.  Then $r'r \cdot x = x$ and $r'r \cdot y = x$. Since $x$, $r\cdot x$, and $y$ are in the same strong orbit, there is some $s\in R_T$ such that $s\cdot (r\cdot x) = y$; then $sr \cdot x = y$ and $sr \cdot y = y$.  Then $r'r\cdot (sr \cdot x) = r'r\cdot y = x$.  On the other hand, $sr \cdot (r'r \cdot x) = sr \cdot x = y$.  Since $(r'r)(sr) = (sr)(r'r)$ it must be the case that $x =y$ and multiplication by $r$ is injective, as desired.

    Next, consider surjectivity.  Let $x\in T$; we will show that $x$ is in the image of $r\cdot$.  Since $x$ and $r\cdot x$ are in the same strong orbit $T$, there exists $r'\in R_T$ such that $x = r'\cdot (r \cdot x)$.  But then 
    \[x = r'r \cdot x = r \cdot (r'\cdot x),\]
    so $r\cdot$ is surjective, as desired.
\end{proof}

\begin{remark}
    Proposition~\ref{lem:injact} is false if $A$ is not abelian. Surjectivity clearly fails, but injectivity fails even in very simple cases.  Consider, for example, the  $\N * \N$-set whose underlying set is $\{0,1,\ldots,10\}$, where the two generators $s$ and $t$ act by $s(i) = \max(i-1,0)$ and $t(i) = \min(i+1,10)$.  These are well-defined actions on an $A$-set with a single non-trivial strong orbit $T$, with $s,t\in R_T$, but the actions by $s$ and by $t$ are not injective.
\end{remark}

%Note that by definition, $A$ acts transitively on any $A$-ambit. \maru{what does this mean here? an $A$-ambit $T$ is not closed under the $A$-action, since a priori we have $aT\subseteq T'$ for some other ambit $T'$} \inna{transitivity generally means that if $x,y\in T$ then there exists $a\in A$ such that $a \cdot x  = y$.  Should I put it in?}\maru{Ok, what confused me is that it feels a little weird to talk about transitivity for sets that are not closed under the action, but you're right that the property is still true. I think we might want to either use a different term, or have a remark about this.}

%\begin{definition}
%    A map $f: S \rto T$ of $A$-sets is a function that commutes with the $A$-action, in the sense that $a \cdot f(x) = f(a \cdot x)$ for all $a\in A$ and $x\in X$.  Such a function is an \emph{isomorphism} if there exists a map $g: T \rto S$ such that $fg = 1_T$ and $gf = 1_S$.  We write $S \cong T$ if there exists an isomorphism from $S$ to $T$.
%\end{definition}

We can now use the above result to describe the set of $A$-automorphisms of any invariant strong orbit, in analogy with the case of transitive $G$-sets.

\begin{theorem}\label{Aut(T)}
  Let $A$ be an abelian monoid, $S$ be an $A$-set, and $T\subseteq S$ be a non-trivial invariant strong orbit. Then $\Aut(T_+) \cong R_T\ss F_T$. 
\end{theorem}

\begin{proof}
Fix any element $x\in T$, and let $r\in R_x=R_T$.  Then $(r\cdot)$ is an automorphism of $T_+$,  as it is a bijection by Proposition \ref{lem:injact} and $A$-equivariant since $A$ is abelian. Moreover, note that if $r\sim_{F_T} r'$ then $(r\cdot) = (r'\cdot)$, and so this automorphism only depends on the equivalence class $[r]\in R_T\ss F_T$.  Now suppose that for $r,s\in R_x$, we have that $(r\cdot) = (s\cdot)$.  Let $y = r\cdot x$, and let $r'\in R_T$ be such that $r' \cdot y = x$.  Then $r'r,r's\in F_x$; in particular, $r(r's) = s (r'r)$, so $r\sim_{F_x} s$.  Thus we have shown that there is an injective function $R_x\ss F_x \rto \Aut(T_+)$.

    It remains to check that this function is surjective.  To show this, we  claim that any automorphism $\varphi:T_+ \rto T_+$ is uniquely determined by $\varphi(x)$. First, suppose that $\varphi(x) = x$.  For any element $y\in T$, there is some $r\in R_x$ such that $r\cdot x = y$, since $x$ and $y$ are in the same strong orbit.  Then $\varphi(y) = \varphi(r\cdot x) = r\cdot \varphi(x) = r\cdot x = y$, so $\varphi$ fixes $y$ as well; we conclude that in this case, $\varphi$ is the identity on $T$.  Now suppose that we have two automorphisms $\varphi,\psi$ such that $\varphi(x) = \psi(x)$.  Then $\psi^{-1}\circ \varphi\in \Aut(T_+)$ fixes $x$, and thus is the identity; in particular, $\varphi = \psi$.  Finally, as any $\varphi(x)\in T$ is in the same strong orbit as $x$, it can be written as  $r\cdot x = \varphi(x)$ for some $r\in R_T$; thus the function $R_x/F_x \rto \Aut(T_+)$ is an isomorphism, as desired.
\end{proof}

\section{Classification of the \texorpdfstring{$K$}{K}-theory of \texorpdfstring{$A$}{A}-sets}

%For $G$,once this classification is known (of the $G$-orbits), one uses the fact that every $G$-set is a disjoint union of $G$-orbits.
In this section, we use our characterization of strong orbits and their automorphisms to give a computation of the $K$-theory of the category of finite $A$-sets for any finite abelian monoid $A$. 

\begin{notation}
    In this section, $A$ is assumed to be finite, and any $A$-set is assumed to also be finite, even if it is not explicitly stated. 
\end{notation}

Following \cite{Hchuck}, we consider the category of finite $A$-sets as an ACGW category for the purposes of $K$-theory. This agrees with its $K$-theory as a Waldhausen or quasi-exact category, but notably, it gives us access to a Devissage theorem that will be crucial for our purposes. Throughout this section, we freely use the definitions and terminology regarding (pre-)ACGW categories from \cite{CZ-cgw}. A reader new to these frameworks might wish to skim through \cite[Section 3.1]{SS-cgw} where the motivating example of exact categories is presented in detail through this double-categorical lens, or \cite{ACGWhomology} which contains a more informal account of ACGW categories through the example of finite sets.

We thus begin by detailing how $A$-sets assemble into an ACGW category. To this end, we remind the reader of an important class of maps between $A$-sets.

\begin{definition}\label{defn.quotient}
    Let $A$ be a monoid and $X\subseteq S$ be $A$-sets. The \emph{quotient} $S/X$ is the pointed set of equivalence classes $S/\sim_X$, where the equivalence relation $\sim_X$ is such that $s\sim_X s$ for all $s\in S$, and $x\sim_X x'$ for all $x,x'\in X$. This is an $A$-set, with action given by $a\cdot [s]\coloneqq [a\cdot s]$. 

    The \emph{quotient map} $S\to S/X$ is the pointed map $s\mapsto [s]$; this is $A$-equivariant by definition of the action on $S/X$. 
\end{definition}

\begin{remark}
    Unlike the quotient construction $\ss$ for submonoids considered in the previous sections, the quotient of $A$-sets is taken in a pointed way. All elements in $X$ are collapsed to the basepoint of $X$, whose equivalence class becomes the basepoint of $S/X$, and all other elements are left untouched; that is, $[s]=\{s\}$ for all $s\in S\backslash X$. 
\end{remark}

\begin{proposition}
 For any abelian monoid $A$, there is an ACGW category  $\ASets$ with 
 \begin{description}
     \item[objects] finite $A$-sets.
     \item[m-morphisms] $A$-equivariant injections.
     \item[e-morphisms] opposites of $A$-equivariant quotient maps. 
     \item[squares] commutative squares.
     \item[distinguished squares] bicartesian squares; i.e.\ squares which are both pushouts and pullbacks. 
     \item[initial object] the trivial $A$-set $\{*\}$.
       \item[kernels] the kernel of an e-morphism $p\colon S\elto S'$ is given by the inclusion $p^{-1}(*)\mrto S$.
     \item[cokernels] the cokernel of an m-morphism $f\colon X\mrto S$ is given by the quotient $S\elto S/ f(X)$.
     \item[restricted pushouts] are pushouts of $A$-sets for m-morphisms, and pullbacks for e-morphisms.
 \end{description}
\end{proposition}

 \begin{proof}
     This proceeds as in \cite[Section 2]{Hchuck}, as the category of finite $A$-sets is also quasi-exact. Note that in \cite{Hchuck} the monoid $A$ is assumed to be Noetherian, as the authors work with the category of finitely generated $A$-sets which are otherwise not closed under subobjects; this hypothesis is not necessary in order to work with finite $A$-sets, and their arguments go through verbatim. 
 \end{proof}

%note for us: we COULD talk about Noetherian A and finitely generated guys, or about any A and "admissible" monos and epis. But who cares right? it doesn't seem worth it since our classification is only for A being finite and for finite A-sets.

\begin{remark}\label{silly.remark}
    The m- and e-morphisms behave as one would expect when it comes to strong orbits, in the sense that an  injection $f\colon X\to S$ will include entire strong orbits of $X$ into entire strong orbits of $S$, and a quotient $p\colon S\to S/X$ will collapse entire strong orbits to the basepoint. 
    
    Indeed, for the first claim, if $x,x'\in X$ are in the same strong orbit, then $x'=r\cdot x$ for some $r\in R_x$ and so $f(x')=f(r\cdot x)=r\cdot f(x)$; since $R_x=R_{f(x)}$ by Lemma \ref{lem:ambittopair}, this shows that $f(x)$ and $f(x')$ are in the same strong orbit in $S$. Conversely, if $f(x)$ is in the same strong orbit as some $s\in S$, then $s=r\cdot f(x)=f(r\cdot x)$ for some $r\in R_{f(x)}=R_x$ and so $s$ is the image of an element $r\cdot x$ in the same strong orbit as $x$. 

    For the second claim, suppose $s,s'\in S$ are in the same strong orbit with $s'=r\cdot s$ for some $r\in R_s$. If $p(s)=[*]$, then $p(s')=p(r\cdot s)=r\cdot p(s)=r\cdot [*]=[*]$ and so $p$ collapses the entire strong orbit of $s$.
\end{remark}

    The $K$-theory of $\ASets$ as an ACGW category agrees with its $K$-theory as a quasi-exact category, as they both admit equivalent $Q$-constructions; see \cite[Definition 4.1]{CZ-cgw} and \cite[Section 3.2]{deitmar06}. In both of these settings, the kernel-cokernel sequences are up to isomorphism all of the form \[X\mrto S\elto S/X\] for an $A$-subset $X\subseteq S$. Notably, these sequences are generally not split,
    as illustrated by the following example.

    \begin{example}
    Consider the monoid $A_d$ from Example \ref{ex:A_d}, and let $S=A_d$ with action given by multiplication, and $X=\{*, t^d,\dots, t^{n-1}\}$. We then get an $A_d$-equivariant injection $X\mrto S$ whose cokernel is $S/X=\{[*],[1],\dots,[t^{d-1}]\}$. However, the quotient map $S\twoheadrightarrow S/X$ does not admit a retraction, as the subset $S/X$ is not an $A_d$-subset of $S$: we have that $t\cdot [t^{d-1}]=[*]$ in $S/X$ but $t\cdot t^{d-1}=t^d\neq *$ in $S$.
    \end{example}

As we have learned in Section \ref{section:classifying.ambits}, we have the most control over the invariant strong orbits of an $A$-set. However, the same behavior that causes cofiber sequences to not be split makes it so that not every $A$-set decomposes as a disjoint union of its invariant strong orbits, as not all strong orbits need be invariant. For instance, in the example above of $S=A_d$, the only invariant strong orbit is the cycle $\{t^d,\dots,t^{n-1}\}$ which evidently does not cover the entire $A$-set $S$.

To circumvent this issue, our strategy will be to consider the full subcategory $\ASets^\inv$ of $\ASets$ of those $A$-sets for which all strong orbits are invariant, and then show that $K(\ASets^\inv)\simeq K(\ASets)$. We start by introducing this subcategory more carefully.

\begin{definition}
    Let $\ASets^\inv$ denote the full double category of $\ASets$ whose objects are the finite $A$-sets $S$ for which every strong orbit is invariant.  
\end{definition}

\begin{remark}
    In other words, the objects of $\ASets^\inv$ are precisely the $A$-sets $S$ such that for all  $a\in A$ and $x\in S$, we have that either $a\cdot x=*$ or that $a\cdot x$ and $x$ are in the same strong orbit. Equivalently, these are the $A$-sets $S$ such that for all $x\in S$ we have that $A\backslash R_x=\{a\in A \mid a\cdot x=*\}$. 
\end{remark}

\begin{proposition}\label{B.is.preACGW}
    The double category $\ASets^\inv$ is a pre-ACGW subcategory of $\ASets$, which is closed under subobjects and quotients, and such that the inclusion $\ASets^\inv\cap \M\to \M$ creates restricted pushouts in the sense of \cite[Definition 6.1]{CZ-cgw}. 
\end{proposition}
\begin{proof}
    As the double category $\ASets^\inv$ clearly includes the initial object $\{*\}$, to prove that it is a pre-ACGW subcategory of $\ASets$ it suffices to show that it is closed under kernels, cokernels, and pullbacks. 

    The claim that $\ASets^\inv$ is closed under kernels means that for any quotient map $p\colon S\to S'$ where both $S,S'\in B$, we have that its kernel $p^{-1}(*)$ is also in $\ASets^\inv$. Dually, the claim that $\ASets^\inv$ is closed under cokernels means that for any injection $f\colon X\to S$ with $X,S\in\ASets^\inv$, we have that its cokernel $S/f(X)$ is in $\ASets^\inv$. It then suffices to show that $\ASets^\inv$ is closed under all subobjects and quotients, as this is a strictly stronger claim; we focus on this instead. 

Note that any object $S\in\ASets^\inv$ must decompose, up to isomorphism, as a finite disjoint union of invariant strong orbits \[S\cong \{*\}\amalg\bigsqcup_{i=1}^n T_i\] for some $n\geq 0$. To show that $\ASets^\inv$ is closed under subobjects, suppose we have an injection $f\colon X\to S$ with $S\in \ASets^\inv$; we want to show that $X\in \ASets^\inv$. By Remark \ref{silly.remark}, the strong orbits in $X$ are, up to isomorphism, a subset of the strong orbits in $S$; hence $X$ must be an object in $\ASets^\inv$. Similarly, to show that $\ASets^\inv$ is closed under quotients, suppose we have a quotient map $p\colon S\to S'$ with $S\in \ASets^\inv$; we want to show that $S'\in \ASets^\inv$. By Remark \ref{silly.remark}, we know that $p$ simply collapses entire strong orbits; hence the strong orbits in $S'$ are, up to isomorphism, a subset of the strong orbits in $S$ and we deduce that $S'\in\ASets^\inv$ as well. 

   % Next, to show that $\B$ is closed under quotients, suppose we have a quotient map $p\colon S\to S'$ with $S\in \B$; we want to show that $S'\in \B$. Without loss of generality, we can assume that $S'=S/X$ for some $A$-subset $X\subseteq S$. Then, given $[s]\in S/X$ and $a\in A$, we want to show that either $a\cdot [s]=[*]$ or $a\cdot [s]$ and $[s]$ are in the same ambit. Since $S\in \B$, we know that either $a\cdot s=*$ (in which case we are done, as $a\cdot [s]=[a\cdot s]=[*]$) or that $a\cdot s$ and $s$ are in the same ambit in $S$. In the latter case, there is some $r\in R_s$ with $a\cdot s=r\cdot s$, and such that there exists some $r'\in R_s$ with $r'r\cdot s=s$. Thus $r'r\cdot [s]=[r'r\cdot s]=[s]$ and we see that $r\in R_{[s]}$; hence $a\cdot [s]=[a\cdot s]=[r\cdot s]=r\cdot [s]$ for some $r\in R_{[s]}$ which means that $a\cdot [s]$ and $[s]$ are in the same ambit in $S/X$. 

  Pullbacks of m-morphisms in $\ASets$ are computed as pullbacks of $A$-sets, so suppose we have a pullback diagram
  \begin{diagram}
      { X \times_Z Y & Y \\ X & Z \\};
      \arrowsquare[right hook->]{}{}{}{}
  \end{diagram}
  where $X,Y,Z\in \B$; we need to show that $X\times_Z Y\in \ASets^\inv$. However, this is immediate since pullbacks preserve injections, thus $X\times_Z Y$ is a subobject of $X\in \ASets^\inv$ and we have proved that $\ASets^\inv$ is closed under subobjects. Dually, pullbacks of e-morphisms are computed as pushouts of the corresponding quotient maps of $A$-sets, and the desired closure is deduced from the fact that $\ASets^\inv$ is closed under quotients. 

  Finally, we must show that the inclusion $\ASets^\inv\cap \M\to \M$ creates restricted pushouts; that is, given a pushout diagram of m-morphisms
  \begin{diagram}
      { X & Y \\ Z & Y\cup_X Z \\};
      \arrowsquare[right hook->]{}{}{f}{g}
  \end{diagram}
  where $X,Y,Z\in \ASets^\inv$, we need to show that $Y\cup_X Z\in \ASets^\inv$. To see this, consider elements $s\in Y\cup_X Z$ and $a\in A$. If $a\cdot s=*$, we are done. Otherwise, we know that $s$ comes from an element of either $Y$ or $Z$ (possibly both); suppose without loss of generality that $s=f(y)$ with $y\in Y$. Since $*\neq a\cdot s=a\cdot f(y)=f(a\cdot y)$, we have that $a\cdot y\neq *$ in $Y$, so as $Y\in \ASets^\inv$ it must be that $a\cdot y$ and $y$ are in the same strong orbit in $Y$. Thus there exists some $r\in R_y$ such that $a\cdot y=r\cdot y$, but then $a\cdot s=f(a\cdot y)=f(r\cdot y)=r\cdot f(y)=r\cdot s$, with $r\in R_y=R_s$ by Lemma \ref{lem:ambittopair} since $f$ is an injection. We conclude that $a\cdot s$ and $s$ are in the same strong orbit in $Y\cup_X Z$, as desired.  
\end{proof}

%\maru{comment for us: initially I tried to prove that $\E\cap\B\to\E$ preserves restricted pushouts (i.e.\ pullbacks of e-mor), so that we could run the usual style of Devissage with a filtration of inclusions, but I actually don't think this is true! A baby case would be to show that if $X,Y\in \B$ then $X\times Y\in\B$, and I think even this is false. If one of you finds a nice counterexample, we could add this as a remark.}

In fact, using similar arguments we can identify another relevant subcategory.

\begin{proposition}\label{prop.single.ambit.type}
    Let $(R,F)$ be a orbital pair in $A$, %$F\leq R\leq A$ be unpointed submonoids as in Theorem \ref{thm.uniqueness}, 
    and let $\ASets^\inv_{R\ss F}$ be the full double subcategory of $\ASets^\inv$ on those objects whose non-trivial strong orbits are all isomorphic to $R\ss F$. Then $\ASets^\inv_{R\ss F}$ is a pre-ACGW subcategory of $\ASets^\inv$.
\end{proposition}

\begin{proof}
    By the same reasoning used in the proof of Proposition \ref{B.is.preACGW}, it suffices to show that the double category $\ASets^\inv_{R\ss F}$ is closed under subobjects and quotients, as it clearly contains the initial object $\{*\}$. Note that any object $S\in\ASets^\inv_{R\ss F}$ must decompose, up to isomorphism, as a finite disjoint union of invariant strong orbits \[S\cong \{*\}\amalg\bigsqcup_{i=1}^n R\ss F\] for some $n\geq 0$. Therefore, using Remark \ref{silly.remark} which ensures that the strong orbits of any subobject  $X\subseteq S$ are a subset of the strong orbits of $S$, and that any quotient map collapses entire strong orbits, we see that $S$ is closed under subobjects and quotients, which concludes the proof. 
\end{proof}

\begin{remark}\label{rmk.trivial.ambit.type}
    The previous result considers non-trivial strong orbit types, but we could analogously consider the full double subcategory $\ASets^\inv_{\{*\}}$ whose objects are finite disjoint unions of invariant singleton strong orbit; this is also a pre-ACGW subcategory of $\ASets^\inv$, by the same argument.
\end{remark}

We now have all the pieces required to show our main result. For convenience, we first introduce the following notation.

\begin{notation}
    For any abelian monoid $A$, let $\Amb_A$ denote the set of isomorphism classes of invariant strong orbits. 
\end{notation}

\begin{remark}\label{another.silly.rmk}
    Note that  $\Amb_A$ is in bijection with the set $OP_A$ of orbital pairs in $A$,% of all pairs $(R,F)$ of unpointed submonoids $F\leq R\leq A$ satisfying the conditions of Theorem \ref{thm.uniqueness}(2), 
    together with the trivial strong orbit. In particular, if $A$ is finite, then so is $\Amb_A$, as there are only finitely many unpointed submonoids of $A$.
\end{remark}

\begin{theorem}\label{BPQish}
   For any finite abelian monoid $A$, we have that 
    \[K(\ASets) \simeq \bigvee_{\substack{T\in \Amb_A\\ T \neq \{*\}}} \Sigma_+^\infty B\Aut(T_+).\]
\end{theorem}

\begin{proof}
 As a first step, we use Devissage (\cite[Theorem 6.2]{CZ-cgw}) to show that \[K(\ASets)\simeq K(\ASets^\inv).\] We have already checked that $\ASets^\inv$ satisfies most of the hypotheses of that theorem in Proposition \ref{B.is.preACGW}, so it remains to show that every finite $A$-set $S\in\ASets$ admits a finite sequence \[S=S_n\elto S_{n-1}\elto S_{n-2}\elto\cdots\elto S_1\elto S_0=\{*\}\] in which the kernel of each quotient map $S_i\elto S_{i-1}$ is an object in $\ASets^\inv$.

If $S$ is the trivial $A$-set $\{*\}$, then we are done, as $\{*\}\in \ASets^\inv$ so we may take the trivial sequence of length 1. Otherwise, let $X_{n-1}\subseteq S$ denote the disjoint union of all the invariant strong orbits in $S$; by \cite[Proposition 2.3(3)]{burnside}, we know that $X_{n-1}$ contains at least one non-trivial invariant strong orbit. %, in addition to the trivial invariant strong orbit $\{*\}$.
Note that $X_{n-1}$ is an $A$-subset of $S$, as it was defined to be closed under the action of $A$, so we may define $S_{n-1}\coloneqq S/X_{n-1}$ and let the map $S\elto S_{n-1}$ be the quotient map. By construction, the kernel of this map is precisely $X_{n-1}$, which is an object of $\ASets^\inv$ as all its strong orbits are invariant.

If $S_{n-1}=\{*\}$, then we have already constructed the desired sequence. Otherwise, we may repeat the process with the $A$-set $S_{n-1}$, letting $X_{n-2}\subseteq S_{n-1}$ be the disjoint union of all the invariant strong orbits in $S_{n-1}$, and letting $S_{n-2}\coloneqq S_{n-1}/X_{n-2}$. The kernel of the quotient map $S_{n-1}\to S_{n-2}$ is then $X_{n-2}\in \ASets^\inv$. Note that in each step of this process, we have that $X_i\neq\{*\}$ by \cite[Proposition 2.3(3)]{burnside}, and so the $A$-set $S_i=S_{i+1}/X_i$ has strictly fewer elements than $S_{i+1}$. Since $S$ is a finite $A$-set, this ensures that the process must terminate after a finite number of steps, producing the desired sequence and proving that $K(\ASets)\simeq K(\ASets^\inv)$.

As a second step, we claim that $\ASets^\inv$ is equivalent to the finite product $$\prod_{T\in\Amb_A} \ASets^\inv_{T},$$ where $\ASets^\inv_T$ denotes the double subcategory $\ASets^\inv_{R_T\ss F_T}$ of Proposition \ref{prop.single.ambit.type} when $T\neq \{*\}$, or the double subcategory $\ASets^\inv_{\{*\}}$ of Remark \ref{rmk.trivial.ambit.type} when $T=\{*\}$. The double functor \[F\colon \prod_{T\in\Amb_A} \ASets^\inv_{T}\to \ASets^\inv\] maps a tuple $(S_T)_{T\in\Amb_A}$, where each $S_T\in\ASets^\inv_T$ is an $A$-set whose only non-trivial strong orbit type is of type $T$, to the $A$-set $\coprod_{T\in\Amb_A} S_T$, which is a disjoint union of the $A$-sets $S_T$ in which all the basepoints have been identified. It is straightforward to see that this is indeed a double functor, that is, it preserves m- and e-morphisms and squares in a functorial way. Moreover, this preserves distinguished squares, is surjective on objects, and fully faithful on m- and e-morphisms and on squares with a fixed boundary. In particular, this implies that $F$ induces an equivalence of categories on the $Q$-construction of both of these pre-ACGW categories, from which we obtain an equivalence 
\[K(\ASets^\inv)\simeq K(\prod_{T\in\Amb_A} \ASets^\inv_{T})\simeq \bigvee_{T\in\Amb_A} K(\ASets^\inv_{T}).\]

The final step is then to prove that $K(\ASets^\inv_T)\simeq \Sigma_+^\infty B\Aut(T_+)$ for each strong orbit type $T\in\Amb_A$.  For this, first note that each $\ASets^\inv_T$ is equivalent to its skeleton $\sk\ASets^\inv_T$, which can be chosen so that each non-trivial object is a finite coproduct of copies of $T_+$, where basepoints are identified; hence $K(\ASets^\inv_T)\simeq K(\sk\ASets^\inv_T)$. Now, $\sk\ASets^\inv_T$ thus defined is a pre-ACGW category in which all m- and e-morphisms are (up to isomorphism) simply coproduct inclusions. By \cite[Propositions 5.9 and 5.12]{SS-cgw}, we get that $K(\sk\ASets^\inv_T)$ as a pre-ACGW category agrees with $K^\otimes ((\sk\ASets^\inv_T)^\iso)$ as a symmetric monoidal category, where  $(\sk\ASets^\inv_T)^\iso$ denotes the maximal groupoid of $\sk\ASets^\inv_T$.

We now claim that the symmetric monoidal category $(\sk\ASets^\inv_T)^\iso$ is equivalent to the symmetric monoidal category $(\GSets^\mathrm{free})^\iso$ of finite free $G$-sets for the finite abelian group $G=\Aut(T_+)$. To see this, we define a functor $H\colon (\sk\ASets^\inv_T)^\iso\to (\GSets^\mathrm{free})^\iso$ as follows. On objects, we have that the objects of the first are finite coproducts of the generator $T_+$, while the objects of the second are all isomorphic to a finite coproduct of the group $G=\Aut(T_+)$; hence the correspondence $H(\coprod_n T_+)=\coprod_n \Aut(T_+)$ is essentially surjective. On morphisms, it suffices to specify the action of $H$ on an isomorphism $\varphi\colon T_+\to T_+$, i.e.\ an element $\varphi\in\Aut(T_+)$; this is mapped to the isomorphism $\Aut(T_+)\to \Aut(T_+)$ given by multiplication by $\varphi$. The verification that this is fully faithful follows standard arguments, analogous to the ones in the proof of Theorem \ref{Aut(T)}. We conclude that \[K^\otimes ((\sk\ASets^\inv_T)^\iso)\simeq K^\otimes((\GSets^\mathrm{free})^\iso),\] where the latter is equivalent to $\Sigma_+^\infty BG=\Sigma_+^\infty B\Aut(T_+)$ by the Barratt--Priddy--Quillen--Segal theorem \cite[Proposition 3.6]{segal_coh}.
\end{proof}

\begin{theorem}\label{mainthm}
For any finite abelian monoid $A$, we have that 
      \[K(\ASets) \simeq \bigvee_{\substack{(R,F)\in OP_A}} \Sigma^\infty_+ B(R\ss F).\] %$F\subseteq R$ range among all unpointed submonoids of $A$ satisfying the conditions of Theorem \ref{thm.uniqueness}.
\end{theorem}
 \begin{proof}
     This is a consequence of Theorem \ref{BPQish}, together with Remark \ref{another.silly.rmk} and Theorem \ref{Aut(T)}. 
 \end{proof}

\section{Functoriality of the splitting}\label{functoriality}

Theorem \ref{mainthm} can be rephrased in the following manner.  There is a CGW functor \[u:\prod_{(R,F)} \ASets^{\inv}_{R \ss F} \rto \ASets\] given by taking a tuple $\{S_{R\ss F}\}_{(R,F)}$ to the $A$-set $\coprod_{(R,F)} S_{R\ss F}$.  Theorem \ref{mainthm} states that $K(u)$ is a weak equivalence.  The question becomes: how functorial is this when $A$ changes?

The first case we consider is the case when $f$ is injective:

%To begin we consider the question of when a monoid homomorphism $f$ induces homomorphisms on orbits.  
%\begin{lemma} \label{lem:ifpairmaps}
%   Let $f:A \rto B$ be a monoid homomorphism. Suppose that $(R,F)\in OP_A$ is such that $(f(R),f(F)) \in OP_B$.  Then $f$ induces a group homomorphism $\tilde f:R \ss F \rto f(R) \ss f(F)$ and therefore $f$ induces a functor 

%%%Maru: I commented out the centered equation below:
%5    \[\ASets^\inv_{R \ss F} \rto \BSets_{f(R)\ss f(F)}^\inv.\]
%\end{lemma}

%\begin{proof}
%    That $\tilde f$ is a group homomorphism follows directly from the fact that $f$ is a monoid homomorphism.  Since $f$ induces a group homomorphism $R \ss F \rto f(R) \ss f(F)$, it follows that any strong orbit that is isomorphic to $R \ss F$ can be mapped to $f(R) \ss f(F)$ by picking $x\in R \ss F$ mapping it arbitrarily, and then using the $A$-action to map the other points in the orbit.  Any two such choices are canonically isomorphic, so this will give a well-defined functor. \inna{CHECKE THIS}  
%\end{proof}

%\maru{I'm a little unsure about how this functor acts on morphisms in the case when $f$ is not injective... Can you describe it? Or, alternatively, I'd be ok with just removing lemma 5.1 and jumping straight to 5.2}

%However, it is not obvious how to ensure that the condition of this lemma holds.  Here is one example in which it does. 

\begin{lemma} \label{lem:goodinj}
Let $f:A \rto B$ be an injective pointed homomorphism of monoids.  Suppose in addition that we can consider $A$ to be a $B$-set by defining $b \cdot a = a'a$ if $b = f(a')$, and $b \cdot a = *$ otherwise.  Then there is a CGW functor $f_*:\ASets \rto \BSets$ given by extending every $A$-action on a set $S$ to a $B$-action in an analogous way. 
\end{lemma}

\begin{proof}
    Let $S$ be an $A$-set.  We claim that we can consider $S$ to be a $B$-set by defining for $b\in B$ and $x\in S$, 
    \[b \cdot x = \begin{cases}
        a \cdot x \caseif b= f(a) \\
        * \caseotherwise.
    \end{cases}\]
    To prove this it suffices to check that if $a\cdot x \neq *$, then $a \neq bb'$ for any $b,b'\notin f(A)$.  Thus suppose for the sake of contradiction that this is the case, so that $a = bb'$ for $b,b'\notin f(A)$.  Then when we consider $A$ as a $B$-set as in the condition of the lemma, we have $a \cdot 1 = a \neq *$, but $b' \cdot (b \cdot 1) = *$, a contradiction.  Thus this is not possible, and the functor is well-defined on objects.  That it is well-defined on morphisms follows immediately.%, since morphisms are simply equivariant functions.
\end{proof}

The fact that this is well-defined on $\ASets$ extends to it being well-defined on invariant strong orbits, and thus implies that the splitting is functorial in the following sense:

\begin{lemma}
    Suppose $f:A \rto B$ satisfies the conditions of Lemma~\ref{lem:goodinj}.  Then there is a commutative diagram
    \begin{diagram}
        {  |[text depth=0.5ex]|\prod_{(R,F)\in OP_A} \ASets^\inv_{R \ss F} & |[text depth=0.5ex]|\prod_{(R,F) \in OP_B} \BSets^\inv_{R \ss F} \\
        \ASets & \BSets \\};
        \we{1-1}{2-1} \we{1-2}{2-2}
        \to{1-1}{1-2}^{\prod f_*} \to{2-1}{2-2}^{f_*}
    \end{diagram}
\end{lemma}

\begin{proof}
It suffices to show that there is a CGW functor $f_*: \ASets_{R\ss F}^\inv \rto \BSets^\inv_{f(R)\ss f(F)}$ for each $(R,F)\in OP_A$.
    The key point is that in this case every orbital pair $(R,F)\in OP_A$ satisfies the condition that $(f(R),f(F))\in OP_B$.  We check the conditions in turn.  By definition $(f(R),f(F))$ are unpointed monoids and $f(F) \leq f(R)$.  To check that $f(R)$ is saturated, suppose that $bb'\in f(R)$.  Write $a = bb'$; since $R$ is unpointed, $a\cdot a \neq *$.  If $b\notin f(A)$ then $b \cdot a = *$ by definition.  But then when we consider $A$ as a $B$-set, $a \cdot a = bb' \cdot a = b' \cdot (b \cdot a) = *$, a contradiction.  Thus we cannot have $b\notin f(A)$, and therefore $f(R)$ is saturated.  Since $f$ is injective, $f(R) \ss f(F) \cong R \ss F$, which is a group.  Lastly, if $b\in B$ and $y\in f(F)$ such that $by\in f(F)$ then $by\in F(R)$, which implies that $b\in f(A)$.  Since $(R,F)$ is a orbital pair, $b\in f(F)$, as desired.
\end{proof}

There is also another form of functoriality that one may want to consider. Let  $f: A \rto B$ be any homomorphism of pointed monoids.  There is an induced CGW functor $f^*: \BSets \rto \ASets$ given by sending a $B$-set $S$ to the $A$-set defined by the action $a \cdot x = f(a) \cdot x$.  This gives the following solid-arrow diagram:
\begin{diagram}
    { |[text depth=0.5ex]| \prod_{(R,F)\in OP_B} \BSets^{\inv}_{R \ss F} &|[text depth=0.5ex]| \prod_{(R,F)\in OP_A} \ASets^{\inv}_{R \ss F} \\
    \BSets & \ASets \\};
    \we{1-1}{2-1} \we{1-2}{2-2} \to{2-1}{2-2}^{f^*}
    \diagArrow{densely dotted,->}{1-1}{1-2}
\end{diagram}
The natural question becomes: does the dotted arrow exist?

\begin{proposition}
    Let $f: A \rto B$ be a surjective map of monoids.  Let $(R,F)\in OP_B$; then $(f^{-1}(B),f^{-1}(A))\in OP_A$.
\end{proposition}

\begin{proof}
    We must check that $f^{-1}(R)$ and $f^{-1}(F)$ are unpointed submonoids with $f^{-1}(F) \leq f^{-1}(R) \leq A$ such that $f^{-1}(R)$ is saturated, $f^{-1}(R) \ss f^{-1}(F)$ is a group, and such that for any $a\in A$ and $y\in f^{-1}(F)$ if $ay \in f^{-1}(F)$ then $a\in f^{-1}(F)$.  We take these in turn.

    That $f^{-1}(R)$ and $f^{-1}(F)$ are unpointed follows directly from the fact that $f$ is a pointed homomorphism.  That they are submonoids follows from the fact that if $x,y\in f^{-1}(R)$ (resp. $f^{-1}(F)$) then $f(x),f(y)\in R$ (resp. $F$) and thus $f(xy) = f(x)f(y)\in R$ (resp. $F$), so that $xy\in f^{-1}(R)$ (resp. $F$).  Suppose that $xy\in f^{-1}(R)$.  Then $f(xy) = f(x)f(y) \in R$ which implies (since $R$ is saturated) that $f(x),f(y)\in R$, which implies that $x,y\in f^{-1}(R)$, as desired.  Thus $f^{-1}(R)$ is saturated.  Now suppose that $x\in f^{-1}(R)$.  To prove that $f^{-1}(R)\ss f^{-1}(F)$ is a group it suffices to show that there exists $y\in f^{-1}(R)$ such that $xy \in f^{-1}(F)$.  Let $y'\in R$ be such that $f(x)y' \in F$.  Since $f$ is surjective, there exists $y\in f^{-1}(y')$; for this $y$, $xy\in f^{-1}(F)$, as desired. Lastly, suppose that $a\in A$ and $y\in f^{-1}(F)$ such that $ay\in f^{-1}(F)$.  Then $f(a)f(y)\in F$ and therefore $f(a)\in F$.  Thus $a\in f^{-1}(F)$, as desired.
\end{proof}

\begin{remark}
  The surjectivity condition is needed to ensure that $f^{-1}(R)\ss f^{-1}(F)$ is a group. Otherwise, given an element $[x]$ in this quotient with $f(x)\in R$, the fact that $R\ss F$ is a group implies that there is an inverse $[r]\in R\ss F$ of $[f(x)]\in R\ss F$, but without surjectivity it may very well be that every representative of the equivalence class of $[r]$ has an empty preimage under $f$.
\end{remark}

%This proposition has an interesting issue: the homomorphism $f^{-1}(R) \ss f^{-1}(F) \rto R \ss F$ goes in the \emph{opposite} direction of $f^*$. Thus in order to construct the dotted map, we must somehow reverse the arrow.

\begin{proposition} \label{transport}
    If $f:A\to B$ is surjective and normal, i.e.\ it is injective away from the basepoint, then the dotted arrow exists.  Moreover, this arrow induces a splitting $$K(\ASets) \simeq K(\BSets) \vee X.$$
\end{proposition}

\begin{proof}
    If $f$ is normal then $f^{-1}(R) \cong R$ and $f^{-1}(F) \cong F$, and in particular %the induced homomorphism $f^{-1}(R) \ss f^{-1}(F) \rto R \ss F$ is an isomorphism.
    we have an isomorphism $f^{-1}(R) \ss f^{-1}(F) \cong R \ss F$. Thus the dotted arrow can be defined to take each $B$-set in the component $\BSets^\inv_{R \ss F}$ to its isomorphic counterpart in the component $\ASets^\inv_{f^{-1}(R)\ss f^{-1}(F)}$.

    Moreover, we see that every component indexed by $(R,F)\in OP_B$ maps isomorphically onto a unique component indexed by $(R,F) \in OP_A$.  This proves that $K(\ASets) \simeq K(\BSets) \vee X$, where $X$ is weakly equivalent to the $K$-theory of those components  that are not indexed by a orbital pair of the form $(f^{-1}(R), f^{-1}(F))$ for some $(R,F)\in OP_B$.
\end{proof}

\section{Applications and Examples}\label{applications}

\subsection{General observations} We record some simple yet interesting observations that can be deduced from the classification of the $K$-theory of $A$-sets.

\begin{proposition} \label{prop:countingspheres}
    For any finite abelian monoid $A$, $K(\ASets) \simeq \mathbb{S}^{\# OP_A} \vee X$, where $\pi_0 X = 0$.
\end{proposition}

\begin{proof}
    For any discrete group $G$, $\Sigma^\infty_+ BG \simeq \mathbb{S} \vee \Sigma^\infty BG$, where $\pi_0\Sigma^\infty BG = 0$.  Thus we can rewrite the formula in Theorem~\ref{mainthm} as
    \[K(\ASets) \simeq \mathbb{S}^{\# OP_A} \vee \bigvee_{(R,P)\in OP_A} \Sigma^\infty B\Aut(R \ss F).\]
\end{proof}

From this we get a minimum bound on the $K$-theory of $A$-sets:

\begin{proposition}\label{always.a.sphere}
 For any finite abelian monoid $A$, we have that  $K(\ASets)$ contains at least one copy of the sphere spectrum, and at least two copies if $A$ is non-trivial and has at least two units, or has no zero divisors.
\end{proposition}
 \begin{proof}
    It suffices to prove this by showing that $OP_A$ is always nonempty, and has at least two elements if the conditions of the lemma hold. 

    First, note that the pair $(A^\times,A^\times)$ is always an orbital pair, so $OP_A$ is always nonempty.  

    Secondly, consider two cases.  If $A^\times \neq \{1\}$ then the pair $(A^\times,\{1\})$ is also an orbital pair; this handles the first case.
    
    If this is not the case but $A$ has no zero divisors, we claim that the pair $(A \setminus \{*\}, A \setminus\{*\})$ is also an orbital pair.  To see that these are submonoids, let $a,b\in A\setminus\{*\}$; then $ab\neq *$ since $A$ has no zero divisors, and we see that $A\setminus\{*\}$ is closed under multiplication. The remaining conditions of an orbital pair are immediate.
 \end{proof}

% \begin{corollary}
%    For any finite abelian monoid $A$, we have that 
%     \[K(\ASets) \simeq \mathds{S}\vee \bigvee \Sigma_+^\infty B(R\ss F),\] where $F\leq R<A$ range among all unpointed submonoids of $A$ satisfying the conditions of Theorem \ref{thm.uniqueness}.
% \end{corollary}
% \begin{proof}
%     This is a consequence of Theorem \ref{BPQish}, together with Remark \ref{another.silly.rmk} and Theorem \ref{Aut(T)}. \maru{$\Sigma_+^\infty B\{*\}$ is the sphere, so the statement above is correct, right? Also I think I used the wrong font for the sphere}
% \end{proof}

\subsection{Finite abelian groups}

When our monoid is of the form $A=G_+$ for a finite abelian group $G$ (following the notation of Example \ref{groups.as.monoids}), our formula recovers the classical computation for $K(\GSets)$. The fact that the $K$-theory of $G$-sets where $G$ is viewed as a group and the $K$-theory of pointed $G_+$-sets where $G_+$ is viewed as a monoid agree is not new, but we record it here nevertheless as a quick sanity check.

\begin{proposition}\label{K(G)}
    For any finite abelian group $G$, we have that $K(\GSets)\simeq K(G_+\text{-}\mathbf{Sets})$.
\end{proposition}
\begin{proof}
    The classical computation expresses the spectrum on the left as \[K(\GSets) \simeq \bigvee_{H \leq G} \Sigma^\infty_+ B(G/H).\] On the other hand, our computation for the spectrum on the right yields
    \[K(G_+\text{-}\mathbf{Sets}) \simeq \bigvee_{\substack{F \subseteq R \subseteq G_+}} \Sigma^\infty_+ B(R\ss F)\] where $(R,F)\in OP_{G_+}.$ %$F\subseteq R$ range among all unpointed submonoids of $G_+$ satisfying the conditions of Theorem \ref{thm.uniqueness}. 
    Now, recall that the submonoids $F,R$ arise as the $x$-reversible elements $R=R_x$ and the stabilizer $F=F_x$ for any element $x$ in any orbit $T$ whose orbit type is represented by $R\ss F$  (see the proof of Theorem \ref{thm.uniqueness} for details). But for the monoid $G_+$, we always have that $R_x=G$, as all elements in $G$ are invertible and hence $x$-reversible. Moreover, an unpointed submonoid $F\subseteq G$ is precisely a subgroup of $G$, and any such subgroup will give a relative pair; %satisfy the required conditions of Theorem \ref{thm.uniqueness}; 
    thus the two formulas above agree.
\end{proof}

\subsection{Monoids generated by a single element}

Let $A$ be a finite monoid with a single generator. It is not hard to convince oneself that $A$ is always abelian, and moreover, that depending on the relations imposed on the generator, $A$ must be isomorphic to one of the following:
\begin{enumerate}
    \item $(\mathbb{Z}/n\mathbb{Z})_+$; that is, $A=\{*,1,t,\dots,t^{n-1} \mid t^n=1\}$;
    \item $\mathbb{N}/t^n\mathbb{N}$; that is, $A=\{*,1,t,\dots,t^{n-1}\mid t^n=*\}$;
    \item $A_d=\{*,1,t,\dots,t^{n-1}\mid t^n=t^d \text{ for some } d\geq 1\}$.
\end{enumerate}

We can use this classification to compute $K(\ASets)$. 

\begin{proposition}
    Let $A$ be a finite monoid with a single generator. Then 
    \[K(\ASets)\simeq\begin{cases} K(\mathbb{Z}/n\mathbb{Z}) \caseif A\cong (\mathbb{Z}/n\mathbb{Z})_+;\\
    \mathds{S} \caseif A\cong \mathbb{N}/t^n\mathbb{N};\\
    \mathds{S}\vee\bigvee_{k\vert (n-d)} \Sigma^\infty_+ B(\mathbb{Z}/k\mathbb{Z}) \caseotherwise\end{cases}\]
\end{proposition}
\begin{proof}
In the first case we have that $A$ is a pointed group, and so its $K$-theory agrees with $K(\mathbb{Z}/n\mathbb{Z})$, as recalled in Proposition \ref{K(G)}. 

In the second case, we have that $A$ is \emph{partially cancellative}: for any elements $a,b,c\in A$, the fact that $ac=bc\neq *$ implies that $a=b$ (\cite[Definition 1.1]{Hchuck}). Moreover, we claim that any finite $A$-set $S$ must be partially cancellative: for any $x\in S$ and $a,b\in A$, the fact that $a\cdot x=b\cdot x\neq *$ implies that $a=b$. 
To see this, suppose that $t^i\cdot x=t^j\cdot x\neq *$ for some $x\in S$ and some $i<j<n$. Then \[t^i\cdot x=t^j\cdot x=t^{j-i}\cdot (t^i\cdot x)=t^{j-i}\cdot (t^j\cdot x)=t^{j-i}t^{j-i}\cdot (t^i\cdot x)=\dots\] and we can continue raising the power of $t^{m(j-i)}$ until eventually we reach $t^n=*$, which contradicts the assumption that $t^i\cdot x\neq *$. Thus, as all $A$-sets are partially cancellative, we have that $K(\ASets)=K'(A)\simeq \mathds{S}$, where the equality is by definition of $K'$ and the equivalence is given in \cite[Lemma 3.3]{Hchuck}.

Lastly, in the case when $A$ is isomorphic to the monoid $A_d$ of Example \ref{ex:A_d}, we use Theorem \ref{mainthm} to obtain  \[K(A_d\text{-}\mathbf{Sets}) \simeq \bigvee_{\substack{F \subseteq R \subseteq A_d}} \Sigma^\infty_+ B(R\ss F)\] where $F\subseteq R$ range among all unpointed submonoids of $A_d$ such that $R$ is saturated, $R\ss F$ is a group, and $F$ is such that whenever $af\in F$ with $a\in A$, $f\in F$, it must be that $a\in F$. 

We first claim that the only saturated unpointed monoids $R\subseteq A_d$ are either $R=\{1\}$ or $R=A_d\setminus\{*\}$. Indeed, suppose that $R$ contains any element $t^i\neq 1$; since $t^i=tt^{i-1}$ and $R$ is saturated, we have that $t\in R$ and then by closure, every power of $t$ is also in $R$. Next, note that the only pair $F\subseteq R$ in which $R=\{1\}$ participates is the trivial pair $\{1\}\subseteq\{1\}$; this gives the sphere $\mathds{S}$ in our statement, as explained in Corollary  \ref{always.a.sphere}. When $R=A_d\setminus\{*\}$, for any $F\subseteq R$ as in the theorem we see that $R\ss F$ is generated by $[t]$. Since $[t]^n=[t]^d$ and $R\ss F$ is a group, it must be $[t]^{n-d}=[1]$ and so the order of $[t]$ divides $n-d$. Conversely, each divisor $k$ of $n-d$ gives rise to a unique unpointed submonoid $F_k=\{t^i \mid i\equiv 0 \ (\text{mod } k)\}$ satisfying the required conditions and with the property that $R\ss F_k$ is a group of order $k$. We then obtain the formula in the statement. 
\end{proof}

\subsection{Nilpotent ideals}
 Recall that a finite monoid $A$ is \emph{nilpotent} if there exists some natural number $n$ such that $a^n=*$ for all $1\neq a\in A$.  

%\maru{this first lemma is in response to a Q by Christian H, who was like, this must be true right?}
 \begin{lemma}\label{nilpotent}
     If $A$ is a finite nilpotent abelian monoid, then $K(\ASets)\cong\mathds{S}$.  
 \end{lemma}
 \begin{proof}
     We know by Proposition \ref{always.a.sphere} that $K(\ASets)$ contains at least one copy of the sphere spectrum, indexed by the pair of unpointed submonoids $R=F=\{1\}$. We claim that these are the only unpointed submonoids of $A$. To see this, let $R\subseteq A$ be any subset closed under multiplication, and suppose there is an element $1\neq r\in R$. Then $*=r^n\in R$, and $R$ must be pointed. We thus conclude our result by appealing to Theorem \ref{mainthm}. 
 \end{proof}

 This result can be extended to consider quotients of $A$ by nilpotent ideals, in analogy with the classical result for Noetherian rings where we have that $G(R)\simeq G(R/I)$ for any nilpotent ideal $I\leq R$; see for instance \cite[Corollary V.4.2]{weibel_kbook}. In our context, recall that a non-empty subset $I\subseteq A$ is an \emph{ideal} if we have that $a\cdot x\in I$ for any $a\in A$, $x\in I$ (in particular, this implies that $I$ must be pointed). Although ideals are typically not submonoids (as they do not contain the unit; otherwise $I=A$), they are $A$-sets with action given by multiplication, and it makes sense to consider the quotient $A/I$ as in Definition \ref{defn.quotient}. This quotient is itself a monoid, and so it makes sense to study the $K$-theory of $A/I$-sets, as we do in the following result. 

 \begin{proposition}\label{nilpotent.quotient}
     Let $A$ be a finite abelian monoid and $I\subseteq A$ be a nilpotent ideal. Then $$K(\ASets)\simeq K(A/I\text{-}\mathbf{Sets}).$$
 \end{proposition}
\begin{proof}
The case when $I=A$ is precisely the content of Lemma \ref{nilpotent}; we henceforth assume that $I\subsetneq A$.

    Since the projection to the quotient $\pi\colon A\to A/I$ is a surjective normal map, we already know by Proposition \ref{transport} that $K(\ASets)\simeq K(A/I\text{-}\mathbf{Sets})\vee X$. Inspecting the proof of this result, the spectrum $X$ is itself a wedge, indexed relative pairs $(R,F)$ %by the pairs of unpointed submonoids $F\subseteq R\subseteq A$ that satisfy  the conditions of Theorem \ref{thm.uniqueness} 
    and do not arise as preimages under $\pi$ of corresponding unpointed submonoids in $A/I$. We aim to show that there are no such additional pairs.

    To prove this, start with a pair of unpointed submonoids $F\subseteq R$ of $A$ with the required conditions. It suffices to show that $F\cap I=R\cap I=\emptyset$, as then we have isomorphisms (of unpointed submonoids) $F\cong \pi(F)$ and $R\cong \pi(R)$ and it is not hard to show that the pair of unpointed submonoids $\pi(F)\subseteq \pi(R)$ in $A/I$ is a relative pair %satisfies all the conditions of Theorem \ref{thm.uniqueness}. 
    Now, suppose there is an element $x\in R\cap I$; since $I$ is nilpotent, we have that $x^n=*$ for some $n$, but $R$ is closed under multiplication and so $*\in R$, contradicting the fact that $R$ is unpointed. The same argument applies to $F$, and concludes our proof.
    \end{proof}

    \subsection{Local commutative algebras}

    Let $\Bbbk$ be a finite field, and $A$ be a finite local commutative algebra over $\Bbbk$. More explicitly, $A$ is a commutative $\Bbbk$-algebra, which is finite-dimensional as a $\Bbbk$-vector space, and which has a unique maximal ideal $\mathfrak{M}\subsetneq A$. In particular, the multiplication in $A$ allows us to consider $A$ as an abelian monoid, where the basepoint is given by $0$. We can then show the following.  This result is in the tradition of the results of \cite[Lemma 3.3]{Hchuck} and \cite[Remark 2]{deitmar06}, showing that in certain cases invariants of monoid sets are isomorphic to those of their groups of units.

    \begin{proposition}\label{prop.comm.local.algebras}
        For any finite local commutative algebra $A$ over $\Bbbk$, we have that $$K(\ASets)\simeq K(A^\times\text{-}\mathbf{Sets})$$ where $A^\times$ denotes the group of units of $A$.
    \end{proposition}
    \begin{proof}
        Since $A^\times$ is a group, we know by Proposition \ref{K(G)} that $K(A^\times\text{-}\mathbf{Sets})\simeq K(A^\times_+\text{-}\mathbf{Sets})$, where the first considers sets with an action by the group $A^\times$, and the second considers pointed sets with an action by the monoid $A^\times_+$. 
        
        Let $\mathfrak{M}$ denote the unique maximal ideal of $A$. If $A/\mathfrak{M}$ denotes the pointed quotient as in Definition \ref{defn.quotient}, which is itself a monoid, we have a monoid isomorphism $A/\mathfrak{M}\cong A^\times_+$. Indeed, any unit $a\in A^\times$ must be in $A\setminus\mathfrak{M}$ (as having a unit in $\mathfrak{M}$ would imply that $\mathfrak{M}=A$), and conversely, any element $a\in A\setminus\mathfrak{M}$ must be a unit, as otherwise we would get a proper ideal $\mathfrak{M}\subsetneq \mathfrak{M}\cup\langle a\rangle\subsetneq A$ which contradicts the maximality of $\mathfrak{M}$. Therefore, we have that $K(A^\times_+\text{-}\mathbf{Sets})\simeq K(A/\mathfrak{M}\text{-}\mathbf{Sets})$.

        Lastly, the fact that $A$ is finite-dimensional over $\Bbbk$ implies that $A$ is an Artinian ring, and hence its unique maximal ideal must be nilpotent (see for instance the discussion preceding Proposition 8.6 in \cite{atiyahMcdonald}). We thus conclude our proof by Proposition \ref{nilpotent.quotient}.
    \end{proof}

\bibliographystyle{alpha}
\bibliography{SSZ}

\end{document}